\documentclass[twoside,12pt]{article}

\usepackage{amsmath}
\usepackage{amssymb}
\usepackage{amscd}
\usepackage{amsthm}

\numberwithin{equation}{subsection}
 
\theoremstyle{plain}

\newtheorem{Th}[equation]{Theorem}
\newtheorem{Prop}[equation]{Proposition}
\newtheorem{Le}[equation]{Lemma}
\newtheorem{Cor}[equation]{Corollary}

\theoremstyle{remark}

\theoremstyle{definition}

\newtheorem{Def}[equation]{Definition}
\newcommand{\E}{\mathbb E}

\newcommand{\R}{\mathbb R}

\def\ga{\gamma}

\def\ra{\rightarrow}

\def\e{\emph}
\def\i{\infty}
\def\p{\partial}
\def\b{\begin}
\newcommand{\ol}{\overline}

\begin{document}

\title{    \flushleft{\bf{Groups Acting on $CAT(0)$ Square Complexes}}      }
\date{  }
\maketitle

\vspace{-5mm}

\noindent 
Xiangdong Xie\newline
Department of Mathematics,
Washington University,
St.Louis, MO 63130.\newline
Email:   {  xxie@math.wustl.edu}



\pagestyle{myheadings}

\markboth{{\upshape Xiangdong Xie}}{{\upshape   Groups Acting on $CAT(0)$ Square Complexes}}

      
\vspace{3mm}

\noindent
{\small {\bf Abstract.}
We study groups   acting on $CAT(0)$ square complexes.  
In particular we show if $Y$
 is a    nonpositively curved (in the sense of A. D. Alexandrov)
finite square complex 
  and  the vertex links of $Y$ contain no  simple loop consisting of five edges,
 then any  subgroup  of $\pi_1(Y)$ 
either is  virtually free abelian  or    contains
a free  group  of rank two.  In addition we  discuss when   a  group
 generated by two hyperbolic isometries  contains a free group of rank two  and when two points in
 the ideal boundary of a $CAT(0)$ $2$-complex   at Tits distance $\pi$ apart are the endpoints of a geodesic
    in  the $2$-complex.}

\vspace{3mm}



\vspace{3mm}

\setcounter{section}{1}
\setcounter{subsection}{0}

\subsection{Introduction}

We study groups acting isometrically and  cellularly on $CAT(0)$ square complexes.
  A $CAT(0)$ square complex is a simply connected 2-complex such that each 2-cell 
 is isometric to the unit square in the  Euclidean  plane and all simple loops in the vertex links contain
 at least four edges. A $CAT(0)$ square complex  with the induced path metric is a $CAT(0)$
 space.  


Let $G$ be a group acting properly and cocompactly on a $CAT(0)$ square complex by
cellular isometries.   By a result of G. Niblo-L. Reeves (\cite{NR}),  $G$  does not have  
Kazhdan's property (T).  Besides this result,  not much  restriction 
 is known about the  properties of $G$.   $G$ can be an 
  incoherent
    group (for instance $F_2\times F_2$, 
\cite{S}, 
\cite{Bi}),
    a  nonresidually finite group (
\cite{W})     or  even   a 
torsion-free simple  group (
\cite{BM}).   In this paper 
     we are    mainly  interested in the structure of subgroups of $G$. Specifically 
 we would like to know whether a subgroup contains a free group of rank two if it is not virtually
 free abelian.

The  same question can be asked for  any group $G$  acting properly and cocompactly 
 on  a  $CAT(0)$ space $X$.
It   has  been  answered   affirmatively  for the following spaces: 
    trees (\cite{PV}) or more generally  Gromov hyperbolic spaces  (\cite{G});
  certain cubical complexes (\cite{BSw}); 
   Euclidean buildings of rank $\ge 3$ or symmetric spaces (\cite{T}); spaces with 
isolated flats (\cite{HR});        
    Hadamard $3$-manifolds and  certain  real analytic Hadamard $4$-manifolds (\cite{X1}).
The question  in the general case  appears hard.  In particular it is still open
  for Hadamard $4$-manifolds  and  $CAT(0)$ $2$-complexes.  It  is   not even known (\cite{Sw}) 
 whether $G$ has  an infinite  subgroup where each element is of finite order. 
S. Adams and W. Ballmann (\cite{AB}) showed any amenable subgroup of the group $G$ is virtually free
 abelian.  When  $X$ is  a piecewise smooth $2$-complex where
 each edge is contained in at least two $2$-cells,  W. Ballmann and M. Brin (\cite{BBr}) showed
 either $G$ contains  a free group of  rank two or $X$ is isometric to the Euclidean plane.

Below we describe the results in this paper. 
  A square complex is a    \e{quarter plane}  if it  is isomorphic to 
 the first quadrant $\{(x,y)\in \R^2: x\ge 0,\; y\ge 0\}$ with the obvious square complex structure.
  The point   corresponding to the 
  origin is   called  the cone point of the quarter plane. 
  A square complex    homeomorphic to the plane 
   is   a  \e{fake plane}  if    it  is  
the union of five quarter planes
   having  the same cone point  and  disjoint  interiors.   
 The main result of this paper can be stated as follows.

\b{Th}\label{main}
{Let $G$ be a group acting properly and cocompactly on a 
   $CAT(0)$ square complex   $X$ by cellular isometries. If there is no fake plane in $X$, then
 any subgroup   $H$ of $G$ either is virtually free abelian or contains a free group of rank two.}
\end{Th}

In general it is not easy to check whether a  $CAT(0)$ square complex contains a fake plane.
But  if a $CAT(0)$ square complex contains a fake plane, then some 
vertex link (at the cone point) must contain a simple loop with 5 edges.  Notice 
 there are only a finite number of vertex links modulo $G$ 
 and each vertex link is a finite graph.

\b{Cor}\label{vertexl}
{Let $G$ be a group acting properly and cocompactly on a 
   $CAT(0)$ square complex   $X$ by cellular isometries.   If no  vertex link   of $X$ 
 contains  simple loops    consisting   of   five edges,   then
 any subgroup of $G$ either is virtually free abelian or contains a free group of rank two.}
\end{Cor}

   If   each      simple loop in   the vertex links   of  $X$ 
         consists    of   an even number of edges  Corollary \ref{vertexl}    
   also   follows from 
 W. Ballmann  and  J. Swiatkowski's result (\cite{BSw}) on foldable cubical complexes. 
    Recall  for  any odd integer $m\ge 5$,  there   are $CAT(0)$ square complexes  admitting proper
  and cocompact cellular isometric actions such that vertex links contain
  simple  loops  consisting of $m$   edges: 
  (
\cite{Be}, 
\cite{BH}, 
 \cite{BBr2})
  given  any  finite  connected graph $L$ where simple loops
   contain at least 4 edges,  there exists a  $CAT(0)$ square complex admitting   a 
  proper and cocompact cellular isometric action      such that all vertex links are isomorphic to $L$.



For any $CAT(0)$ space $X$,  $d_T$ denotes the 
  Tits metric on the ideal boundary $\p_\i X$.
  For each geodesic $c: R\ra X$, the positive and negative directions of $c$ give rise to
 two points  in the ideal
 boundary, called the endpoints of the geodesic (see Section \ref{idealb}).    When  $g$ 
 is a hyperbolic isometry of a $CAT(0)$ space, we denote  the endpoints of the axes of $g$ 
   by   $g(+\i)$  and  $g(-\i)$.

In order to prove Theorem \ref{main}, we need to decide when a subgroup of $G$ contains a free group
 of rank two.  By W. Ballmann's theorem (see  Section \ref{rank1} or  \cite{B}), 
   if a subgroup $H<G$ contains a rank one isometry
 then $H$ either is virtually infinite cyclic or contains a free group of rank two. 
K. Ruane (\cite{R}) showed if $g, h$ are two hyperbolic isometries and if 
  $d_T(\xi,\eta)>\pi$   for    any $\xi\in \{g(+\i), g(-\i)\}$    and any 
$\eta\in \{h(+\i), h(-\i)\}$, then 
the group generated by  $g$ and $h$ contains a free group of rank two.  Notice   if 
 $\xi, \eta \in \p_\i X$  with $d_T(\xi,\eta)>\pi$,   
  then there is a geodesic in $X$ with $\xi$  and  $\eta$ 
 as endpoints. 
Thus   the following  result should be considered as a generalization of K. Ruane's theorem
 in the case of square complexes. 

\vspace{3mm}

\noindent
{\bf{Theorem \ref{criterion}.}}
\e{Let $X$  be a locally compact $CAT(0)$ square complex.  Suppose $g$ and $h$ are two hyperbolic isometries 
 of $X$ such that  for any $\xi\in \{g(+\i), g(-\i)\}$ and 
any $\eta\in \{h(+\i), h(-\i)\}$, there is a geodesic in $X$   with  $\xi$ and 
 $\eta$   as endpoints.  
 Then  the group   generated by  $g$  and  $h$   contains  a free group 
 of rank two.}

\vspace{3mm}

  R. Alperin,  B. Farb   and  G. Noskov (\cite{AFN})   provided sufficient  conditions 
  for   two hyperbolic isometries to generate a free group of rank two. 
 Let $g,h$ be two hyperbolic isometries of a $CAT(0)$ space  $X$ and $a_g$, $a_h$ their axes.
  In the case the two axes are disjoint  and if $p\in a_g$, $q\in a_h$ realize the distance between 
  the axes, their condition requires     the four angles that   the geodesic segment $pq$ makes  
 with the axes  to    be  $\pi$.      A similar  idea was also used by  W. Ballmann  and  M. Brin
in an earlier paper \cite{BBr}.  Their conditions  are   too stringent for our purposes. Our condition 
is weaker and  more applicable, but only in the square complex case.  It  is  an open  question whether 
Theorem \ref{criterion}  holds for  all  $CAT(0)$ spaces. 


To apply Theorem \ref{criterion} we need to know when  
two points  $\xi, \eta\in \p_\i X$ 
    are the endpoints of a geodesic in $X$.
Recall a necessary condition is   $d_T(\xi, \eta)\ge \pi$,  
 and  a sufficient condition is     $d_T(\xi, \eta)> \pi$. 
 We provide a criterion  for $\xi$    and  $\eta$ 
to be the endpoints of a geodesic in $X$ when
    $d_T(\xi, \eta)=\pi$
and $X$  is a $CAT(0)$ 2-complex. 

A point in the ideal  boundary of a $CAT(0)$ space is called a \e{terminal point}
 if it does not lie in the interior of any Tits geodesic.

\vspace{3mm}

\noindent
{\bf{Theorem \ref{pivisi}.}}
\e{Let   $X$ be a    $CAT(0)$  2-complex    that  
 admits a  proper, cocompact   action by cellular isometries.
If $\xi, \eta\in \p_\i X$    are not terminal points   and  $d_T(\xi, \eta)\ge \pi$,
   then there is a geodesic
 in $X$  with    $\xi$ and $\eta$  as endpoints.}

\vspace{3mm}

The conclusion  of Theorem \ref{pivisi}  does not hold  if $X$ is not a $CAT(0)$ $2$-complex.  For instance, 
 the universal covers of    nonpositively  curved  
 $3$-dimensional graph manifolds (\cite{BS}, \cite{CK}) 
 are counterexamples.

 To complete the proof of Theorem \ref{main}, we need to find two hyperbolic isometries $g, h$
 of the subgroup $H$ such that   $d_T(\xi, \eta)\ge \pi$  for 
 any $\xi\in \{g(+\i), g(-\i)\}$ and 
any $\eta\in \{h(+\i), h(-\i)\}$.
This is done with the help of the following result, whose proof
 is essentially due to W. Ballmann and S. Buyalo (\cite{BB}).  Recall 
for any group $H$ acting isometrically on a $CAT(0)$ space $X$, 
the limit set $\Lambda(H)\subset \p_\i X$ (see Section \ref{action} for definition)  is 
  closed and invariant under $H$. A \e{minimal set } is a closed and $H$-invariant subset 
 of $\Lambda(H)$ that does not properly  contain any closed and $H$-invariant subset 
 of $\Lambda(H)$.  Minimal set always exists when the $CAT(0)$ space is locally compact.

\vspace{3mm}

\noindent
{\bf{Theorem \ref{dis}.}}(\cite{X1})
\e{Let $X$ be a    locally compact $CAT(0)$  space and $H$ a group acting on $X$ by isometries.
  If   $H$   does not contain any  rank one isometry,   then 
 for any minimal set $M$  we have  $d_T(m, \xi)\le \pi$  
    for any  $m\in M$ and any $\xi\in \Lambda(H)$.}


\vspace{3mm}

The paper is organized as follows.  In Section \ref{pre} we recall basic facts about $CAT(0)$ square 
 complexes  
  and collect  results that shall be needed later on.  The topics covered in this section include:
   $CAT(0)$ $2$-complexes, 
  Tits boundary, action of a hyperbolic isometry on the ideal boundary, convex subgroups and minimal set.
In Section \ref{freesg} we prove Theorem \ref{criterion}.  In Section \ref{visibility}
we recall results from \cite{X2} about Tits boundary of $CAT(0)$  $2$-complexes and prove 
Theorem \ref{pivisi}.  In Section \ref{sectionfake}
 we use support sets to show  a  $CAT(0)$ square complex  contains a fake plane if and only if
 the Tits boundary contains   a  simple loop with length $2.5\pi$.
The proof of Theorem \ref{main} is completed in Sections \ref{gelement}
 and \ref{secsingular}.

\noindent
\textbf{Acknowledgment.} \textit{I would 
   like to thank Bruce Kleiner,    Quo-Shin Chi  and   Blake Thornton for helpful 
discussions and suggestions.}

\subsection{Preliminaries}\label{pre}

In  this section   we  first recall basic definitions concerning  $CAT(0)$ spaces  
  and  then record some   results that shall be needed later on.
 The reader  is referred to \cite{B}, \cite{BH}, \cite{R}, \cite{Sw}  and \cite{X1}   
    for more details on the material in this section.

\subsubsection{$CAT(0)$   Square Complexes} \label {cat0square}

Let $Z$ be an arbitrary metric space. For any subset $A\subset Z$ and any $\epsilon>0$, 
 the $\epsilon$-neighborhood of $A$ is  
$N_\epsilon(A)=\{z\in Z: d(z, a)\le \epsilon\; \text{for some } a\in A\}$.
  For any $z\in Z$ and any $r>0$,  
   $\ol{B}(z,r)=\{z'\in Z: d(z,z')\le r\}$  is 
the closed  metric ball with center $z$ and radius $r$.
  For a  group  $G$ of isometries of $Z$,  the action of $G$ on $Z$  is  \e{proper}  if for   any 
   compact set $K\subset Z$,   the set $\{g\in G: g(K)\cap K\not=\phi\}$
  is finite.

   A   \e{minimal geodesic} in $Z$ is  an isometric embedding 
  $c:I\ra Z$ where $I\subset R$ is an interval in the real line $R$.  A \e{geodesic} in $Z$
 is a locally isometric embedding $c:I\ra Z$, that is, each point $t\in I$ has a neighborhood
  $U$ in $I$ such that $c_{|U}$  is a minimal geodesic.  If $c: [a, b]\ra Z$ is a geodesic, 
  then we  say  $c$ is  a \e{geodesic segment}  from $c(a)$ to  $c(b)$; sometimes we also say 
 $c$ connects the two points $c(a)$ and $c(b)$.   We also  abuse language and call
 the image of a geodesic a geodesic. 

 The \e{Euclidean cone} over $Z$ is the
 metric space   $C(Z)$   defined as follows.
 As a set $C(Z)=Z\times [0,\i)/{Z\times \{0\}}$.  The image of $(z,t)$ is denoted   by   $tz$.
$d(t_1z_1, t_2z_2)=t_1+t_2$ if $d(z_1, z_2)\ge \pi$,  and 
$d(t_1z_1, t_2z_2)=\sqrt{t_1^2+t_2^2-2t_1t_2\cos(d(z_1,z_2))}$ if $d(z_1, z_2)\le \pi$.
    The Euclidean  cone  $C_r(Z)$  over $Z$ with radius $r>0$ is the 
 metric ball $\ol{B}(O,r)\subset C(Z)$ where $O=Z\times \{0\}$ is the cone point
  of $C(Z)$.

A $2$-dimensional $CW$-complex   is called a \e{polygonal complex}  if\newline
\noindent
(1) all the attaching maps  are homeomorphisms;\newline
\noindent
(2) the intersection of any two closed cells  is either empty or exactly
 one closed cell.

A $0$-cell  is also called a vertex.   A   closed $2$-cell whose boundary contains $n$ vertices
 is called a $n$-gon.

A  polygonal complex is   \e{piecewise Euclidean}  if there is a metric on each closed cell such that\newline
\noindent
(1) each closed cell is isometric to a  closed convex subset of the Euclidean plane;\newline
\noindent
(2) if $A$, $B$ are two closed cells with $A\subset B$, then the inclusion of $A$ into $B$ 
 is an isometric embedding.

It follows that a $n$-gon in a piecewise Euclidean  polygonal complex   $X$ is isometric to a  convex 
 Euclidean $n$-gon, although the interior angles at some vertices could be $\pi$. 
  We assume there are only a finite number of isometry types for the closed cells of 
$X$ and  always equip  $X$ with the induced path metric.  It follows from \cite{BH}
  that $X$ is a complete metric space. 

Let $X$ be a locally finite  piecewise Euclidean  polygonal complex 
and $x\in X$.   The link 
$Link(X, x)$  is a  metric graph  defined as follows.
    Let $A$ be a closed $2$-cell containing $x$. The  unit tangent space $S_x A$  of $A$ at $x$  is 
isometric to the unit circle with length $2\pi$.  We first define a subset 
  $Link(A,x)$  of $S_x A$. 
 For any $v\in S_x A$, 
  $v\in Link(A,x)$   if and only if 
the initial segment of the geodesic with initial point $x$ and initial direction $v$ lies in
   $A$.  Thus   $Link(A, x)=S_x A$  if $x$ lies in the interior of $A$;  $Link(A,x)$ is a 
  closed semicircle (with length $\pi$)
 if $x$ lies in the interior of a $1$-cell contained in $A$; and 
   $  Link(A,x)$  is a closed segment with length $\alpha$
 if $x$ is a vertex of $A$ and the interior angle 
 of $A$ at $x$ is $\alpha$.   Similarly  if $x$ is contained in a closed $1$-cell $B$ we can define 
  $S_x B$ and $Link(B, x)\subset S_x B$. We note $S_x B$ consists of two points at distance $\pi$
 apart, $Link(B, x)=S_x B$  if $x$ lies in    the  interior of $B$ and $Link(B, x)$ consists of a single point
 if $x$ is a vertex of $B$.  When $x$   lies  in a closed $1$-cell $B$ and $B$ is contained in a 
   closed $2$-cell $A$,  $S_x B$ and $Link(B, x)$ can be naturally identified with subsets of 
  $S_x A$ and $Link(A, x)$ respectively.

We define 
$Link(X, x)=\cup_A Link(A,x)$, where $A$ varies over all closed  $1$-cells and $2$-cells containing $x$. 
 Here $Link(B, x)$ is identified with a subset of $Link(A, x)$ as indicated in the last paragraph when
$x$   lies  in a closed $1$-cell $B$ and $B$ is contained in a 
   closed $2$-cell $A$.
  We let $d_x$ be   the induced path metric on $Link(X,x)$. 


Since $X$ is piecewise Euclidean and locally finite,
  it   is   not hard to see that  for each $x\in X$
 there is some $r>0$ such that $\ol{B}(x, r)\subset X$ is isometric to the 
  Euclidean cone   $C_r(Link(X,x))$   over $Link(X,x)$ with radius $r$.

We say $X$ is \e{nonpositively curved}
 if for each vertex $v\in X$, $Link(X, v)$  contains no simple loop with length strictly less than
  $2\pi$.   A \e{$CAT(0)$ $2$-complex}  in this paper shall     always  
    mean a simply connected,   
  nonpositively curved   and  
 locally finite  piecewise Euclidean  polygonal complex.


 A \e{$CAT(0)$ square complex} is  a $CAT(0)$ $2$-complex  such that  each closed $2$-cell is isometric to
the square  in the   Euclidean  plane with edge length 1.  Let $X$ be a
   $CAT(0)$ square complex. Then for each vertex $v\in X$,  
   each edge in $Link(X, v)$   has length $\pi/2$. Thus each simple loop in $Link(X, v)$    
    contains at least $4$ edges.



   $CAT(0)$ $2$-complexes  are examples of $CAT(0)$ spaces (see for example \cite{BH}). 
Let $X$ be a $CAT(0)$ space.
 One key feature
 of a $CAT(0)$ space is the  convexity of distance function: if  $c_1: I_1\ra X$, $c_2: I_2\ra X$ 
 are geodesics, then $d(c_1(t_1), c_2(t_2))$ is a convex function defined on $I_1\times I_2$.
   It follows that all geodesics  in $X$ are minimal geodesics,  and 
 there is exactly   one geodesic segment connecting two given points in $X$. 
For any $x,y\in X$,  $xy$   shall denote the unique
geodesic segment connecting $x$ and $y$.

 Let $X$ be a $CAT(0)$ space.  For any three points $x,y,z\in X$  with $y,z\not=x$, 
let $x',y',z'\in \E^2$
 be three points on the  Euclidean  plane    such that  $d(x,y)=d(x',y')$, $d(y,z)=d(y',z')$ and
 $d(z,x)=d(z',x')$. Set $\tilde\angle_x(y,z)=\angle_{x'}(y',z')$, the angle 
 between the two  segments  $x'y'$ and $x'z'$ in the   Euclidean  plane. 
 Let $c_1$ and $c_2$ be geodesic segments from $x$ to $y$ and from $x$ to 
    $z$ respectively.  The convexity
 of distance function implies  $\tilde \angle_x(c_1(t), c_2(t))$ is an increasing function 
of $t$.  Thus we can define the \e{angle} at $x$ between $y$ and $z$ as follows:
 $$\angle_x(y, z)=\lim_{t\ra 0}\tilde \angle_x(c_1(t), c_2(t)).$$
  The angle at $x$ satisfies the triangle inequality:  $\angle_x(y, z)\le  \angle_x(y, w)+\angle_x(w, z)$  for all  $y, z,w\not=x$. And we have $\angle_x(y, z)+\angle_y(z, x)+\angle_z(y, x)\le \pi$
for any 
  three distinct points $x,y,z\in X$.

   Let $X$ be  a $CAT(0)$ $2$-complex.   For any $x,y \in X$ with $x\not=y$, 
 the tangent vector of $xy$ at $x$   is  denoted  by  
$\log_x(y)\in Link(X,x)$. 
     Notice  for    any  $x,y,z\in X$  with $y,z\not=x$,   we have 
    $\angle_x(y,z)=\min\{\pi, d_x(\log_x(y), \log_x(z))\}$.       

Let $X$ be a $CAT(0)$ space.   A  subset  $A\subset X$   is   a \e{convex} subset if  $xy\subset A$ 
 for any    $x, y\in A$.  Let $A\subset X$ be   a
 closed convex subset. 
 The orthogonal projection  onto $A$,  $\pi_A: X\ra A$  can be defined as follows:
  for any $x\in X$ the inequality $d(x, \pi_A(x))\le d(x, a)$  holds  for all 
  $a\in A$.  It follows that for any $x\notin A$ and any $a\not=\pi_A(x)$ we have
  $\angle_{\pi_A(x)}(x, a)\ge \pi/2$. 

\subsubsection{Ideal Boundary of a $CAT(0)$ Space} \label{idealb}

   Let $X$ be a CAT(0)
space, for example a $CAT(0)$   $2$-complex.   
A \emph{ray} starting from $p\in X$ is a geodesic   $\alpha: [0,\infty) \rightarrow X$ with 
$\alpha(0)=p$.   Two rays $\alpha_1$ and $\alpha_2$ are \emph{asymptotic} if  
$d(\alpha_1(t), \alpha_2(t))$  is a bounded function   on the interval  $[0,\infty)$.    The 
\emph{ideal boundary} of $X$ 
 is  the set $\partial X$ of asymptotic  classes of rays in $X$.  
   Set $\ol X=X\cup \partial X$.
 For any $p\in X$ 
and any  $\xi\in \partial X$, there is a unique ray   $\ga_{p\xi}: [0,\i)\ra X$  
  that starts from $p$ and belongs to
$\xi$ (the image of $\ga_{p\xi}$ is  denoted    by  $p\xi$).
 Thus   for any $p\in X$ we can identify $\partial X$ with the set   of rays 
starting from $p$.    Fix  $p\in X$.  Let  $\alpha$, $\alpha_i$ ($i=1, 2, \cdots $)
  be rays  starting from $p$.  
We say    $\{\alpha_i\}_{i=1}^\i$ 
 converges to  $\alpha$   if  $\alpha_i$ converges to $\alpha $ uniformly
on compact subsets of $[0, \infty)$.   Similarly for $x_i\in X$ ($i=1,2,\cdots$)  we say 
 $\{x_i\}_{i=1}^\i$  converges to $\xi\in \ol{X}$  if 
$px_i$ converges to the geodesic segment or ray $p\xi$ uniformly on compact subsets.
 In this way  we define a topology on  $\ol{X}$. 
 It  is    easy to check that this topology is independent of the point $p\in X$ and the induced
  topology on $X$ coincides with the  original  metric topology on $X$.    Both this topology 
 and the induced  topology  on $\partial X$   are  called the \emph{cone topology}. 
  $\partial X$ together with the cone topology 
  is called the \emph{geometric boundary} of $X$, and denoted   by  $\partial_{\infty} X$.


Let   $p\in X$    and   $\xi,\eta\in \partial X$,    it follows from  the
convexity of distance function   that 
  $\tilde \angle_p(\ga_{p\xi}(t), \ga_{p\eta}(t))$  is an increasing function of $t$ and   therefore 
the     following  limit exists:
$\lim_{t\ra \i}\tilde \angle_p(\ga_{p\xi}(t), \ga_{p\eta}(t))$. 
 This limit   is independent of the point $p$  and defined to be the 
Tits angle        $\angle_T(\xi, \eta)$    between $\xi$ and $\eta$.
  The \emph{Tits metric} $d_T$   on  $\p X$  is  the path metric induced by
$\angle_T$.  Notice $d_T$ takes value in $[0,\i)\cup \{+\i\}$. 
 In particular,  by Proposition \ref{Titsb} (iv)  below,   
$d_T(\xi, \eta)=\infty$  if  and only if  $\xi, \eta\in \p X$  are in 
 different   path  components   with respect   to  $\angle_T$.  
    $\p X$  equipped  with the  Tits 
metric $d_T$ is denoted   by  $\partial_T X$. 
The Tits topology  and the cone topology are generally quite different.


Below we collect some basic facts concerning  the Tits metric. 
 For more details please see \cite{BGS}  and \cite{BH}.
For any geodesic $c: R\rightarrow X$ in a $CAT(0)$ space, we   call  the two points 
in $\partial_{\infty} X$ determined by the two rays $c_{|{[0, +\infty)}}$ and  $c_{|{(-\infty,0]}}$
 the endpoints of $c$, and denote them by $c(+\infty)$ and $c(-\infty)$  respectively.

\begin{Prop}\label{Titsb}
{Let $X$ be a locally compact  $CAT(0)$  space,    $\partial_T X$ its Tits boundary,
    and $\xi, \eta,  \xi_i,  \eta_i \in \partial_T X$ \e{($i=1,2,\cdots$)}. \newline
\emph{(i)} $\partial_T X$  is  a  CAT(1) space; in particular, there is no  geodesic loop 
   in    $\p_T X$   with length
   strictly  less than $2\pi$;\newline
\emph{(ii)} If  $d_T(\xi, \eta)>\pi$, then there is a geodesic in 
$X$ with $\xi, \eta$ as endpoints;\newline
\e{(iii)} If $\xi_i\rightarrow \xi$ and $\eta_i\rightarrow \eta$ in the cone topology,
then $d_T(\xi,\eta)\le {\liminf}_{i\rightarrow \infty} d_T(\xi_i,\eta_i)$;\newline
\emph{(iv)} If   $d_T(\xi, \eta)<\infty$,   then there
  is a minimal geodesic  in   $\partial_T X$  from $\xi$ to $\eta$;\newline
\e{(v)} If $d_T(\xi, \eta)\le \pi$,
  then $d_T(\xi, \eta)=\angle_T(\xi,\eta)$.}
\end{Prop}

\subsubsection{Dynamics of Hyperbolic Isometries } \label{dynamics}

Let  $X$ be a $CAT(0)$  space and $g: X\rightarrow X$ an isometry of $X$. 
$g$ is called a \emph{hyperbolic isometry} if it translates a geodesic, that is, if there is a geodesic 
$c: R\rightarrow X$ and a positive number $l$ so that 
$g( c(t))=c(t+l)$  for all $t\in R$;  
  the geodesic $c$ is called an \emph{axis} of $g$.
    All the axes of $g$ are parallel, thus  it  makes sense to denote
    $g(+\infty)=c(+\infty)$,  $g(-\infty)=c(-\infty)$.
Recall two geodesics $c_1, c_2:  R\ra X$ are parallel
 if $d(c_1(t), c_2(t))$ is a bounded function over $R$.

  For a hyperbolic   isometry  $g$, let $Min(g)$ be the union of all the  axes
 of $g$ and $P_g$  the union of all the geodesics parallel 
 to  the  axes of $g$.   
  We call $P_g$  the \e{parallel set} of $g$.   
  Then  $Min(g)\subset P_g$. 
Both   $Min(g)$    and    $P_g$
are   closed  convex  subsets     in  $X$ and 
     split  isometrically as  follows:   
  $Min(g)=Y\times R\subset Y_0\times R=P_g$,  where 
   $Y\subset Y_0$ is a closed    convex   subset of $Y_0$,  each
 $\{y\}\times R$, $y\in Y$  is an axis of $g$ and each $\{y_0\}\times R$, $y_0\in Y_0$
  is parallel to the axes of $g$.  
The geometric boundaries  $\partial_\infty Min(g)$ and $\p_\i P_g$  naturally embed into 
$\partial_\infty X$. 

   Recall each  isometry $g$ of $X$ induces a homeomorphism (with respect to the cone topology) of   $\ol{X}$, 
   which   we   still  denote   by   $g$. 

\begin{Th}{\emph{(V. Schroeder \cite{BGS},  K. Ruane  \cite{R})}}\label{hyperbolicf}
{Let $X$ be a    locally compact $CAT(0)$  space and $g$ a hyperbolic isometry of $X$.  Then:\newline
 \e{(i)} The fixed point set of $g$ in $\partial_\infty X$   is 
  $\partial_\infty Min(g)$;\newline
\e{(ii)} If  $\xi\in \p_\i X-\p_\i P_g$   and  $\eta$ is an accumulation point
 of the set $\{g^i(\xi): i\ge 1\}$,  then  $\eta\in \p_\i P_g$ 
  and $\angle_T(\eta, g(-\i))=\angle_T(\xi, g(-\i))$.}

\end{Th}

\subsubsection{Convex Subgroups} \label{convexsg}

The following definition is due to E. Swenson (\cite{Sw}).

\begin{Def}\label{dcon}
{Let  $X$ be a $CAT(0)$  space and $G$ a group acting  properly    and cocompactly by
 isometries on $X$. A subgroup $H<G$ is a \e{convex subgroup} if there is a closed 
 convex subset $Y\subset X$ such that $h(Y)=Y$ for any $h\in H$ and  $Y/H$ is compact.}
\end{Def}

\begin{Th}{\emph{(E. Swenson \cite{Sw})}}\label{swenson}
{Let $X$ be a $CAT(0)$ space and $G$ a group acting  properly    and cocompactly by
 isometries on $X$.  Suppose $H, K\subset G$ are  convex subgroups  with  $A,B\subset X$ 
the corresponding closed convex subsets.  Then  $H\cap K$ is a convex subgroup. 
Furthermore, 
for any $\epsilon>0$ with  $N_\epsilon(A)\cap  N_\epsilon(B)\not=\phi$, 
$N_\epsilon(A)\cap  N_\epsilon(B)$ is a closed convex  subset corresponding to 
 $H\cap K$.}

\end{Th}

  Recall for any group $G$ and any $g\in G$, 
the centralizer  of $g$ in $G$ is   defined as follows:
 $C_g(G)=\{\gamma\in G: \gamma g=g\gamma\}$.

\begin{Th}{\emph{(K. Ruane  \cite{R})}}\label{mincentral}
{Let $X$ be a $CAT(0)$ space,  $G$ a group acting    properly and cocompactly by
 isometries on $X$ and $g\in G$ a hyperbolic isometry.  Then  $C_g(G)$
is a convex subgroup with 
$Min(g)$   the corresponding closed convex subset.}

\end{Th}

\subsubsection{Action on the Limit Set} \label{action}

Let $X$ be a $CAT(0)$ space and $H$  a group acting on $X$ by isometries.
 A point 
$\xi\in \partial_{\infty} X$ is a \emph{limit point}  of $H$ if 
there is a sequence of elements  $\{h_i\}_{i=1}^{\infty}\subset H$
   such  that  $h_i(x) \rightarrow \xi$ 
for some (hence any) $x\in X$. 
 The \emph{limit set} $\Lambda(H)\subset \partial_{\infty} X $ 
of $H$ is the set of limit points of  $H$. 
  It    is   easy to check that  
$\Lambda(H)$ is closed (in the cone topology) and $H$-invariant.

\begin{Def}\label{minimal}
{A     nonempty closed and $H$-invariant subset $M\subset \Lambda(H)$ is  \emph{minimal }
if it does not contain any proper subset that is  nonempty, closed and $H$-invariant.}
\end{Def}

When the $CAT(0)$ space $X$ is locally compact, the geometric boundary $\partial_{\infty}X$
    and all its closed subsets   are  compact.  It follows from Zorn's lemma that 
   minimal set always exists when $X$ is locally compact.     

 The   proof  of the following result is 
essentially due to  Ballmann and Buyalo (\cite{BB}), 
  who   assumed     $\Lambda(H)= \partial_{\infty} X $.  See Section \ref{rank1}
   for the definition of a rank one isometry.

\begin{Th}{\emph{(\cite{X1})}}\label{dis}
{Let $X$ be a    locally compact $CAT(0)$  space and $H$ a group acting on $X$ by isometries.
  If   $H$   does not contain any  rank one isometry,   then
$d_T(m, \xi)\le \pi$ for any   minimal  set $M$   and any 
$m\in M$, $\xi\in \Lambda(H)$.}
\end{Th}

\subsection{Free Subgroups}\label{freesg}
    
The goal of this section is to establish  a 
criterion (Theorem \ref{criterion}) for the existence of free subgroups
  in a group acting isometrically 
on a $CAT(0)$ square complex.

\subsubsection{Rank One Isometries} \label{rank1}

Let  $X$ be a $CAT(0)$ space. 
      A \emph{flat half-plane}  in $X$  is 
   the  image of an isometric embedding  
  $f:\{(x,y)\in \E^2: y\ge 0\}\rightarrow X$, and in this case we say the geodesic 
  $c:R\rightarrow X$,  
$c(t)=f(t,0)$   bounds the 
flat half-plane.   

\begin{Def}\label{rank1iso}
{A hyperbolic isometry $g$ of  a $CAT(0)$ space $X$ is  called  a \emph{rank  one isometry} 
if   no  axis of $g$  bounds  a  flat half-plane in  $X$.}
\end{Def}


\begin{Th}{\emph{(W. Ballmann  \cite{B})}}\label{dyn}
{Let $X$ be a locally compact $CAT(0)$  space and 
  $g$      a rank one isometry of  $X$. 
   Given any neighborhoods
$U$ of $g(+\infty)$ and $V$ of $g(-\infty)$ in $\ol  X$,  there is an $n\ge 0$ 
such that  $g^k(\ol  X- V)\subset U$ and
$g^{-k}(\ol  X- U)\subset V$  whenever
$k\ge n$.}
\end{Th}

 Theorem \ref{dyn} implies  $g(+\infty)$ and $g(-\infty)$  are the only fixed points 
 of   a rank one isometry $g$ in $\ol  X$. 
   The theorem  also    has   the following  two corollaries.

\begin{Cor}\label{faf}
{Let $X$ be  a locally compact $CAT(0)$  space, 
  $G$     a  group of isometries of  $X$    and 
   $g\in G$    a rank one isometry.  
     Then one of the following holds:\newline
\emph{(i)}  $G$ has a fixed point in $\p_\i X$;\newline
\emph{(ii)}  Some axis $c$ 
  of $g$  is  $G$-invariant; \newline 
\emph{(iii)}   $G$ contains  a free group of rank two.}
\end{Cor}

\begin{Cor}\label{1possibility}
{Let $G $   be   a  group    acting     properly   and cocompactly by isometries on a $CAT(0)$ space,  and $H$ 
a  subgroup of $G$.    If $H$ contains a rank one isometry, then $H$ either
  is virtually infinite cyclic or contains a free group of rank two.}
\end{Cor}



\subsubsection{Ping-Pong Lemma} \label{ppl}

For any group  $G$ and   any  $g_1, g_2\in G$,  $<g_1, g_2>$  denotes the group generated by $g_1$  and $g_2$. 
We recall the following well-known lemma.

\begin{Le}\label{pingpong}
{Let $G$ be a group acting on  a  set $X$, and  $g_1$, $g_2$   two  elements of $G$. 
If $X_1$, $X_2$ are disjoint subsets of $X$ and for all $n\not=0$, $i\not=j$,
  $g_i^n(X_j)\subset X_i$, then the subgroup $<g_1, g_2>$ is free of rank two.}

\end{Le} 

We will apply the Ping-Pong Lemma in the following setting. 
 Let $X$ be a $CAT(0)$ space and $g_1$, $g_2$ two hyperbolic isometries of $X$.
Then there are geodesics $c_1:R\ra X$, $c_2: R\ra X$ and  numbers
 $a, b>0$   with    $g_1(c_1(t))=c_1(t+a)$ and  $g_2(c_2(t))=c_2(t+b)$ for all $t\in R$.
Let $\pi_1: X\ra c_1(R)$ and $\pi_2: X\ra c_2(R)$ be orthogonal projections 
 onto the geodesics $c_1$ and $c_2$ respectively.   Set $X_1=\pi_1^{-1}(c_1((-\i,0]\cup [a, \i)))$ and
$X_2=\pi_2^{-1}(c_2((-\i,0]\cup [b, \i)))$.   The  following lemma is clear. 

\begin{Le}\label{proj}
{Let $X$, $g_1$, $g_2$, $X_1$ and $X_2$ be as above.  If $X_1\cap X_2=\phi$, then 
 the conditions in the Ping-Pong Lemma are satisfied. In particular, 
  $<g_1, g_2>$ is free of rank two.}

\end{Le}

\subsubsection{Free Subgroup Criterion} \label{fsc}


The following is the main   result   of  Section  \ref{freesg}. 

\begin{Th}\label{criterion}
{Let $X$  be a    locally compact
   $CAT(0)$ square complex.  Suppose $g_1$ and $g_2$ are two hyperbolic isometries 
 of $X$ such that  for any $\xi\in \{g_1(+\i), g_1(-\i)\}$ and 
any $\eta\in \{g_2(+\i), g_2(-\i)\}$, there is a geodesic in $X$ with $\xi$   and 
  $\eta$   as endpoints.  
 Then   $<g_1, g_2>$   
contains  a free group 
 of rank two.}
\end{Th}

Let $c:I\ra X$ be a geodesic defined on an interval I. Then for any $p$ in the interior of $c(I)$, 
   the $d_p$ distance between  the two directions of $c$ 
at $p$ is  at least $\pi$.  Recall $d_p$ is a path metric defined on the link  $Link(X,p)$,       and 
    in general  the two directions of $c$ 
at $p$  may have $d_p$  distance  larger than $\pi$. 

\begin{Def}\label{rgeo}
{A geodesic $c:I\ra X$ is an   \e{R-geodesic} if for each point $p$ in 
the interior of $c(I)$,  the $d_p$  distance   between 
the two directions of $c$ 
at $p$  is $\pi$.}

\end{Def}

\begin{Le}\label{nonrgeo}
{With the assumptions   of   Theorem \ref{criterion}.  If at least one of $g_1$, $g_2$ has 
 an axis that is not an   R-geodesic, then   
   $<g_1, g_2>$ 
 contains  a free group 
 of rank two.}
\end{Le}
\b{proof} Notice   if  a geodesic  bounds a flat half-plane then it 
is an   $R$-geodesic.
 The assumption in the lemma implies at least one of $g_1$, $g_2$ is a rank one isometry.
Now the lemma follows   from   Corollary \ref{faf}.

\end{proof}

  In light of  Lemma \ref{nonrgeo}, we will assume  the axes of $g_1$ and $g_2$ are 
$R$-geodesics.

\b{Le}\label{perp}
{Given any two R-geodesic rays  $c_1:[0, \i)\ra X$, $c_2:[0,\i)\ra X$, there is a  number $a>0$ with the following 
 property:   if $q_i$ \e{($i=1,2$)} is a point in the interior of  $c_i$, and $p\in X$, $p\not=q_1, q_2$
 such that \newline
\e{(i)} $pq_1\cap pq_2=\{p\}$;\newline
\e{(ii)} $pq_1$ and $pq_2$ are R-geodesics;\newline
\e{(iii)} for each $i=1,2$, at least one of the angles $\angle_{q_i}(p, c_i(0))$,    
$\angle_{q_i}(p, c_i(\i))$  is $\pi/2$;\newline
 then $\angle_p(q_1, q_2)\ge a$.}
\end{Le}
\b{proof}
We first notice that for each $R$-geodesic $c$, there is a number $\alpha$, $0\le \alpha\le \pi/4$  
  with the following property:
   if $x,y\in c$ with $x\not=y$ and $xy$ 
  lies in some square $S$ of $X$, then
  $d_x(\log_x(y), \xi)$  takes value in
$\{\alpha, \pi/2-\alpha, \pi/2+\alpha, \pi-\alpha\}$   for  any $\xi$   in  
    $ Link(X,x)$  that is 
    parallel to one of the edges of $S$. 
 Let $\alpha_i$ ($i=1,2$) be such a number corresponding to $c_i$.   Since $pq_i$ is also 
an  $R$-geodesic and at least one of the angles $\angle_{q_i}(p, c_i(0))$,    
$\angle_{q_i}(p, c_i(\i))$  is $\pi/2$,   the number $\alpha_i$ also corresponds to 
  $pq_i$.  The fact  $pq_1\cap pq_2=\{p\}$  implies $\angle_p(q_1, q_2)\not=0$. 
Now it   is  easy to see that 
 $\angle_p(q_1, q_2)\ge |\alpha_1-\alpha_2|$  when $\alpha_1\not=\alpha_2$;
$\angle_p(q_1, q_2)\ge \min\{2\alpha_1, \pi/2-2\alpha_1\}$   when 
 $0<\alpha_1=\alpha_2<\pi/4$;   
 and $\angle_p(q_1, q_2)\ge \pi/2$   when 
$\alpha_1=\alpha_2=0$ or $\pi/4$.   

\end{proof}

\b{Le}\label{smallangle}
{Let $\xi_1, \xi_2\in \p_\i X$    and 
$c_1, c_2: [0, \i)\ra X$    be   two  
     rays belonging to $\xi_1$  and $\xi_2$  respectively.
    Suppose 
there is a geodesic $c$  in $X$  with  $\xi_1$   and  $\xi_2$ as endpoints.  
Then for $i\not=j$,  $\angle_{c_i(t)}(c_i(0),  c_j(t))\ra 0$ as $t\ra \i$.}
\end{Le}
\b{proof}
Since $c$ is a  geodesic  in $X$  with  $\xi_1$    and  $\xi_2$  as endpoints,  there is a number $\epsilon>0$ such that
 $c_i\subset N_\epsilon(c)$ ($i=1,2$).   The convexity of distance function implies 
     for large enough $t$, there is a 
point $p_t\in c_1(t)c_2(t)$ with 
$d(c(0), p_t)\le \epsilon$.  Now the lemma follows by considering the triangle
 $\Delta(c_i(t)c_i(0)p_t)$. 

\end{proof}

 Let $g_1$, $g_2$ be as in Theorem \ref{criterion} and $c_1$, $c_2$  
axes of $g_1$, $g_2$ respectively. By Lemma \ref{nonrgeo}, we may  assume 
     $c_1$, $c_2$   are $R$-geodesics. 
Let $\pi_i:X\ra c_i(R)$ ($i=1,2$) be the orthogonal 
projection onto $c_i$.
 If there is     some   $T>0$ such that  
$$\pi_1^{-1}(c_1((-\i, -T]\cup [T, \i)))\cap \pi_2^{-1}(c_2((-\i, -T]\cup [T, \i)))=\phi,$$ 
  then for large enough $n$, $g_1^n$ and $g_2^n$ satisfy the condition in Lemma \ref{proj}
  and  therefore generate a free group of rank two.
We shall prove  there is some $T>0$ such that $\pi_1^{-1}(c_1([T, \i)))\cap \pi_2^{-1}(c_2([T, \i)))=\phi$, the other three cases are similar. 

\b{Prop}\label{emptyi}
{Let $\xi_1, \xi_2\in \p_\i X$    and 
   $c_1, c_2: [0, \i)\ra X$  be  two   rays belonging to $\xi_1$  and $\xi_2$  respectively.
 Suppose 
  $c_1$ and $c_2$ are $R$-geodesics and 
there is a geodesic $c$ in $X$ with   $\xi_1$   and $\xi_2$  as endpoints.   Then  
$\pi_1^{-1}(c_1([T, \i)))\cap \pi_2^{-1}(c_2([T, \i)))=\phi$  
for large enough $T>0$, 
where 
 $\pi_i: X\ra c_i([0,\i))$ is the orthogonal projection onto $c_i$.}
\end{Prop}
\b{proof}
Since  $c_1$ and $c_2$ are $R$-geodesics, there is a  number $a>0$ with the 
 property  stated in Lemma \ref{perp}.  By Lemma \ref{smallangle}, 
 there is  $T>0$   with  $\angle_{c_1(t)}(c_1(0),  c_2(t))<a/3$  and 
$\angle_{c_2(t)}(c_2(0),  c_1(t))<a/3$  for all $t\ge T$.

Suppose  there  exists  $t>T$ with
$\pi_1^{-1}(c_1([t, \i)))\cap \pi_2^{-1}(c_2([t, \i)))\not=\phi$. We will derive a contradiction from this. 
  Pick
$p\in \pi_1^{-1}(c_1([t, \i)))\cap \pi_2^{-1}(c_2([t, \i)))$ and set $q_i=\pi_i(p)$ ($i=1,2$). 
Let $\sigma_2: [0, b_2]\ra X$ be the geodesic    segment   from $p$ to $q_2$.
   Consider  ${\pi_1}_{|{\sigma_2}}: \sigma_2\ra c_1$ and 
let $t_2=\max\{t\in [0, b_2]:\pi_1(\sigma_2(t))=q_1\}$.  Then the geodesic $q_1\sigma_2(t_2)$
  is an   $R$-geodesic and at least one of the angles 
 $\angle_{q_1}(c_1(0), \sigma_2(t_2))$, $\angle_{q_1}(c_1(\i), \sigma_2(t_2))$ is $\pi/2$.
   By replacing $p$ with $\sigma_2(t_2)$, we may assume  $pq_1$ is an   $R$-geodesic and 
 at least one of the angles $\angle_{q_1}(c_1(0), p)$, $\angle_{q_1}(c_1(\i), p)$ is $\pi/2$.
  Similarly we may assume 
$pq_2$ is an   $R$-geodesic and 
 at least one of the angles $\angle_{q_2}(c_2(0), p)$, $\angle_{q_2}(c_2(\i), p)$ is $\pi/2$.
In general  $pq_1$ and $pq_2$ share an initial segment: $pq_1\cap pq_2=pp'$. By replacing 
  $p$ with $p'$ if necessary we may assume $pq_1\cap pq_2=\{p\}$.  Now the conditions in 
Lemma \ref{perp} are satisfied and  so $\angle_p(q_1, q_2)\ge a$.

  Consider the triangle $\Delta(pq_1q_2)$.  Since $q_1=\pi_1(p)$,   
$\angle_{q_1}(p, c_1(0))\ge \pi/2$. It follows that 
$\angle_{q_1}(p, q_2)\ge \angle_{q_1}(p, c_1(0))-\angle_{q_1}(q_2, c_1(0))\ge \pi/2-a/3$.
  Similarly $\angle_{q_2}(p, q_1)\ge  \pi/2-a/3$.
  Now $\angle_{q_1}(p, q_2)+ \angle_{q_2}(p, q_1)+\angle_p(q_1, q_2)\ge 2(\pi/2-a/3)+a>\pi$,
   a contradiction.

\end{proof}

The proof of Theorem \ref{criterion}
 is now complete.

\subsection{ $\pi$-Visibility}\label{visibility}

Let $X$ be a locally compact $CAT(0)$ space and  $\xi, \eta\in \p_T X$. 
If there is a geodesic in $X$ with $\xi$  and  $\eta$  as endpoints, then  $d_T(\xi,\eta)\ge \pi$. 
 For the converse, Proposition \ref{Titsb}
  says   if  $d_T(\xi,\eta)>\pi$ then there is a geodesic in $X$ with  $\xi$   and  $\eta$
   as endpoints. 
 Our main objective in this section (Theorem \ref{pivisi}) is to provide a sufficient
 condition for the existence of a   geodesic with 
$\xi$   and  $\eta$  as endpoints   when 
  $X$ is a     $CAT(0)$     $2$-complex    and 
$d_T(\xi,\eta)=\pi$.
 We shall also collect some facts about the Tits boundary of $CAT(0)$ square complexes.

\subsubsection{Tits Boundary of $CAT(0)$ $2$-Complexes}  \label{tbo2c}

In this section we recall some facts about the Tits boundary of $CAT(0)$  $2$-complexes. 
The definition of a $CAT(0)$  $2$-complex  is given in Section \ref{cat0square}.
A sector in the Euclidean plane is a closed convex subset whose boundary is the union of
 two rays starting from the origin.  A \e{flat sector}  $S$ in a $CAT(0)$ space  $X$ is the image of 
 an isometric embedding from  a  sector in the Euclidean plane into   $X$. 
 The inclusion  $S\subset X$ clearly induces an isometric embedding  
  $\p_T S\subset \p_T X$.  
The following is a theorem from \cite{X2} adapted to our setting, 
 which   says that away from the endpoints,   a   segment   in the Tits boundary 
  is   the   ideal boundary  of a  flat sector.

\begin{Th}{\e{(\cite{X2})}}\label{t13}   
{Let   $X$ be a   $CAT(0)$  2-complex    that  
 admits a  proper, cocompact   action by cellular isometries.
 Suppose  
$\gamma:[0,h]\rightarrow
\partial_T X$ is  a  geodesic    segment in the Tits 
boundary of $X$ with length 
$h\le \pi$. Then for any $\epsilon > 0$, there exists 
a    flat sector $S$ in $X$    with 
$\partial_T S=\gamma([\epsilon, h-\epsilon])$.}
 
\end{Th}

We also need the following result concerning rays and flat sectors in a $CAT(0)$ 2-complex.

\begin{Prop} {\e{(\cite{X2})}}\label{rays}
{Let   $X$ be a  $CAT(0)$  2-complex     that  
 admits a  proper, cocompact   action by cellular isometries. Suppose $S\subset X$ is a flat sector 
  and $c\subset  X$   is   a     ray.
If $S\cap c=\phi$ and $c$ represents an interior point    of  $\partial_T S$,  then there is a ray 
$c'\subset S\cap X^{(1)}$ asymptotic to $c$,  where    $X^{(1)}$  is the 1-skeleton of $X$. }
\end{Prop}

Let $B$ be the set  $\{(x,y)\in \E^2: x\ge 0, 0\le y\le 1\}$ with 
the obvious square complex structure.  A \e{flat half-strip}  in a $CAT(0)$ square complex  $X$
 is the image of a cellular isometric embedding from $B$ into $X$.

\b{Def}\label{special}
{Let $X$  be a  $CAT(0)$ square complex.  A point $\xi\in \p_T X$ is a \e{singular point}
  if  there is a  flat half-strip in $X$ such that the rays contained in the  flat half-strip 
    belong  to $\xi$.
A  point $\xi\in \p_T X$ is called a \e{regular point} if it  is not a singular point but
  lies in the interior of some geodesic segment in $\p_T X$.}
\end{Def}

\b{Le}\label{sexistence}
{Let $X$  be a  $CAT(0)$ square complex  that  
 admits a  proper, cocompact   action by cellular isometries, 
   and $\sigma\subset \p_T X$ a  geodesic  with length strictly greater than $\pi/2$.  Then there
   exists  a singular point in the   interior of $\sigma$.}
\end{Le}
\b{proof}  We may assume $\pi/2<length(\sigma)<\pi$.  By Theorem \ref{t13} there 
is a flat sector $S$ such that  $\p_T S\subset \sigma$ has length strictly  greater than 
 $\pi/2$. Since $S$ (away from the boundary) is tessellated by squares, it   is   easy to see that 
 there is a singular point in the interior of $\p_T S$.

\end{proof}

\b{Le}\label{bounded}
{Let $X$  be a  $CAT(0)$ square complex   that  
 admits a  proper, cocompact   action by cellular isometries,    and  $\xi, \eta\in \p_T X$ 
 with $d_T(\xi,\eta)\le \pi/2$.  Then there is at most one singular 
point in the interior of  $\xi\eta$.}
\end{Le}
\b{proof}
Suppose there are two singular points $\xi_1$ and $\xi_2$ in the 
interior of  $\xi\eta$.    Theorem \ref{t13} implies there exists a flat sector $S$ such that 
 $\xi_1$ and $\xi_2$ lie in the 
interior of      $\p_T S$   and  $length(\p_T S)<\pi/2$.  
     Since $\xi_1$ and $\xi_2$  are singular points,
  there are two  flat half-strips $B_1$ and $B_2$  and two rays $c_1\subset B_1$,
 $c_2\subset B_2$ asymptotic to $\xi_1$ and $\xi_2$   respectively.   
 If $c_1\cap S\not=\phi$, then a subray of $c_1$ is contained in $S$ and it follows 
$\xi_1$ is represented by a ray contained in $S\cap X^{(1)}$, where 
$X^{(1)}$  is the 1-skeleton of $X$. If $c_1\cap S=\phi$,  then 
 Proposition \ref{rays} implies also that   
$\xi_1$ is represented by a ray contained in $S\cap X^{(1)}$. 
  So in any case $\xi_1$ is represented by a ray contained in $S\cap X^{(1)}$. 
Similarly $\xi_2$ is represented by a ray contained in $S\cap X^{(1)}$. 
We clearly have a contradiction here since $\xi_1\not=\xi_2$ and the length of 
$\p_T S$  is strictly less than $\pi/2$. 

\end{proof}

 Recall (see \cite{X2}) if  $X$ is  a    CAT(0) 2-complex   that   admits an isometric cocompact   action,  
 then for any  $\xi\in \partial_T X  $ and any $r: 0 < r < \frac{\pi}{2}$,
      the   metric  ball  $\ol{B}(\xi, r)$ is an $R$-tree.  It is necessary to
 distinguish those points in  $\p_T X$ that are \lq\lq dead ends''.

\b{Def}\label{terminalp}
{Let $X$ be a $CAT(0)$ 2-complex  and $\xi\in \p_T X$.  We
     say $\xi$ is a \e{terminal point}
 if $\xi$ does not lie in the interior of any geodesic segment in $\p_T X$.}
\end{Def}
 
Lemmas \ref{sexistence}
     and \ref{bounded}  easily imply the following proposition.

\b{Prop}\label{specialig}
{Let $X$  be a  $CAT(0)$ square complex  that  
 admits a  proper, cocompact   action by cellular isometries, 
and 
 $\sigma:[0, l]\ra \p_T X$     be   a  geodesic   segment  with length $l>\pi/2$. 
        If $\sigma(0)$ and $\sigma(l)$  are not terminal points, 
  then there  are numbers  $0\le t_1<t_2<\cdots <t_n\le l$   with the following 
  properties:\newline
\e{(i)} $t_1<\pi/2$, $l-t_n<\pi/2$ and $t_{i+1}=t_i+\pi/2$  for $1\le i\le n-1$;\newline 
\e{(ii)} $\{\sigma(t_1), \cdots, \sigma(t_n)\}$
  is the set of   singular points on $\sigma$.\newline
In particular, if 
$\xi, \eta\in \p_T X$   are  singular   but not terminal points 
    and $d_T(\xi,\eta)< \i$, then $d_T(\xi, \eta)$ is an integral multiple of 
  $\pi/2$.}
\end{Prop}

Notice Proposition \ref{specialig} implies the length of any simple closed
geodesic in $\p_T X$   is  an integral multiple of $\pi/2$.

  Since  small metric balls in $\p_T X$ are $R$-trees,  it makes sense 
to talk about branch points in  $\p_T X$.  By using  Theorem \ref{t13}
 it  is   not hard to prove  the following result:

\b{Prop}\e{(\cite{X2})}\label{branch}
{Let $X$  be a  $CAT(0)$ square complex  that  
 admits a  proper, cocompact   action by cellular isometries.  Then any branch point in 
$\p_T X$  is a singular point.}
\end{Prop}

\subsubsection{$\pi$-Visibility}  \label{vis}

Let   $X$ be  a locally compact $CAT(0)$ 
space and $\xi, \eta\in \p_T X$ with $d_T(\xi,\eta)=\pi$.  
If there is a geodesic  $c$ in $X$   with $\xi$    and  $\eta$   as   endpoints,   then $c$ bounds a flat 
 half-plane (see \cite{BH}).  
In general,    there is  no geodesic in $X$ with 
$\xi$    and  $\eta$   as   endpoints.

\b{Th}\label{pivisi}
{Let   $X$ be a  $CAT(0)$  2-complex    that  
 admits a  proper, cocompact   action by cellular isometries.
If $\xi,\eta\in \p_T X$ are not terminal points and $d_T(\xi,\eta)\ge \pi$, then there is a geodesic
 in $X$   with $\xi$    and  $\eta$   as   endpoints.}
\end{Th}

\b{proof}   By Proposition \ref{Titsb} 
    we may assume $d_T(\xi,\eta)=\pi$.  Let $\sigma:[0,\pi]\ra \p_T X$ be a 
 minimal geodesic from $\xi$ to $\eta$.  Since $\xi,\eta\in \p_T X$ are not terminal points 
 and small metric balls in   $\p_T X$  are $R$-trees, there is some $\epsilon$, $0<\epsilon<\pi/4$ such that 
  $\sigma$  extends to 
 a locally isometric map  $[-\epsilon, \pi+\epsilon]\ra \p_T X$, which is still denoted  by  $\sigma$.

Notice $\sigma_{|[-\epsilon, \pi/2+\epsilon]}$ is a minimal geodesic   in  $\p_T X$
with length less than 
 $\pi$.     By Theorem \ref{t13} there is  a flat sector $S_1$ in $X$ such 
that $\p_T {S_1}=\sigma([-\epsilon/2, \pi/2+\epsilon/2])$.    Similarly there is a 
 flat sector $S_2$ in $X$ with $\p_T {S_2}=\sigma([\pi/2-\epsilon/2, \pi+\epsilon/2])$. 
 Since $X$ is a 2-complex and   
$\p_T {S_1}\cap \p_T {S_2}=\sigma([\pi/2-\epsilon/2, \pi/2+\epsilon/2])$
    is a
 nontrivial interval,   we have  $S_1\cap S_2\not=\phi$.  Pick  $x\in S_1\cap S_2$ and let
 $S\subset S_1\cap S_2$ be the subsector with cone point $x$ and 
 $\p_T S=\sigma([\pi/2-\epsilon/2, \pi/2+\epsilon/2])$.  Fix a point $p$ in the 
 interior of the    flat  sector $S$, and let $c_1$, $c_2$, $c_3$ be rays starting from $p$
   belonging  to  $\sigma(0)$, $\sigma(\pi)$  and  $\sigma(\pi/2)$ respectively.
Since $c_1$ and $c_3$ are contained in the  flat sector $S_1$, $\angle_p(c_1(\i), c_3(\i))=\pi/2$.
 Similarly  $\angle_p(c_2(\i), c_3(\i))=\pi/2$. Since the initial 
segment of each $c_i$ is contained in   the flat sector $S$, it    is   clear that the angle
 $\angle_p(c_1(\i), c_2(\i))=\pi$.  It follows that $c_1\cup c_2$ is a 
complete geodesic  in $X$  with endpoints  $\xi$    and  $\eta$.

\end{proof}

\b{remark}\label{a4coun}
\e{The conclusion  of Theorem \ref{pivisi}  does not hold  if $X$ is not a $CAT(0)$ $2$-complex.  For instance, 
 the universal covers of   nonpositively   curved 
  $3$-dimensional graph manifolds (\cite{BS}, \cite{CK}) 
 are counterexamples.}
\end{remark}





\subsection{Fake Planes and Short  Loops  in the Tits Boundary} \label{sectionfake}

The goal of this section is to prove  Proposition \ref{profake}, which provides a criterion 
 for the existence of a simple geodesic loop
 with length $\frac{5\pi}{2}$ in   the Tits boundary  of a $CAT(0)$ square complex.
 This result shall be used in Sections \ref{gelement}
   and \ref{secsingular}.

\subsubsection{Support Sets of Short Closed Geodesics  in $\p_T X$} \label{support}

Let $X$ be  a   $CAT(0)$  $2$-complex and $x\in X$. There is a map
 $\log_x: \p_T X\ra Link(X, x)$  that sends $\xi\in \p_T X$ to the   tangent   vector  
of 
 ${x\xi}$ at $x$.  Clearly  $\log_x$ is  distance nonincreasing:
   $d_x(\log_x(\xi), \log_x(\eta))\le d_T(\xi,\eta)$ for any $\xi,\eta\in \p_T X$.

\b{Def}\label{supportset} 
{Let $X$ be  a  $CAT(0)$  $2$-complex and $c\subset \p_T X$ a  simple closed geodesic in 
$\p_T X$. The \e{support set}  $supp(c)\subset  X$ of $c$ is   defined  as follows:
  $x\in supp(c)$  if  
the loop $\log_x(c)$   represents a nontrivial  class in $H_1(Link(X,x))$.}
\end{Def}



 Recall  $C_r(c)$ denotes  the Euclidean cone over $c$ with radius $r>0$ (Section \ref{cat0square}).
  Since $c$ is   a  circle,  $C_r(c)$   is  homeomorphic to the  closed unit disk
  in the plane.    Let   $\p C_r(c)$   be the boundary circle of  $C_r(c)$. 
For any $x\in X$ and any $r>0$,  the  map $f_{x,r}: C_r(c)\ra X$  is defined by 
$f_{x,r}(t\xi)=\gamma_{x\xi}(t)$, $\xi\in c$, $0\le t\le r$.       
      Clearly      $f_{x,r}$     represents a class in
 $H_2(X, X-\{x\})$.  By using homotopy along geodesic segments it   is    easy 
   to see that     for any $r, r'>0$, 
  the two maps $f_{x,r}$  and $f_{x,r'}$   represent the same class in 
    $H_2(X, X-\{x\})$,  that is,   $[f_{x,r}]=[f_{x,r'}]\in H_2(X, X-\{x\})$.

\b{Le}\label{equi}
{Let $X$ be  a  locally compact   $CAT(0)$  $2$-complex, $x\in X$  and $c\subset \p_T X$ a  simple closed geodesic in 
$\p_T X$.   Then   $x\in supp(c)$   if   and  only  if 
     $f_{x,r}$  represents a nontrivial class in
 $H_2(X, X-\{x\})$   for  any $r>0$.  }
\end{Le}
\b{proof}
Since $X$ is   locally compact  and piecewise Euclidean, 
for each $x\in X$ and each $r>0$ there is some $r'$   with $0<r'\le r$ such that 
 the closed metric ball $\ol{B}(x,r')$ 
 is isometric
 to the Euclidean cone $C_{r'}(Link(X, x))$ over $Link(X,x)$ with radius $r'$. 
   By  the remark before the lemma,   $[f_{x,r}]=[f_{x,r'}]\in H_2(X, X-\{x\})$.


Set $B=\ol{B}(x,r')$. 
The exact sequence of the pair $(B, B-\{x\})$
 gives us an 
isomorphism $\p:H_2(B, B-\{x\})\ra H_1(B-\{x\})$. 
Note $H_1(B-\{x\})\cong H_1(\p B)$ and the map 
$\log_x: \p B\ra Link(X,x)$ which sends $y$ to the direction of $xy$ at $x$ is 
 a homeomorphism.
Now observe that $\p$ sends 
$[f_{x,r'}]\in H_2(X, X-\{x\})\cong H_2(B, B-\{x\})$ to 
$[\log_x(c)]\in H_1(Link(X,x))\cong H_1(\p B)
\cong H_1(B-\{x\})  $.

\end{proof}

\b{Le}\label{simpleloop}
{If a simple closed geodesic $c$ in $\p_T X$ has length strictly less than $4\pi$  and
 $x\in supp(c)$,  then $\log_x(c)$ is homotopic to a simple closed geodesic 
 in $Link(X,x)$.}
\end{Le}

\b{proof}
Note $Link(X,x)$ is a finite metric graph where each simple loop has length
 at least $2\pi$.
The loop $\log_x(c)$ is homotopic to a    closed geodesic  $c_x$  in  $Link(X, x)$  whose length does not
 exceed the length of $\log_x(c)$.   If $c_x$ is not simple, then its length
 must be at least $4\pi$. The lemma follows.

\end{proof}

 Let $c$ be  a  simple closed geodesic in $\p_T X$ with 
 length strictly  less than  $4\pi$   and  $x\in supp(c)$.  
By Lemma \ref{simpleloop} 
$\log_x(c)$  is homotopic to a simple loop $c_x$ in  $Link(X,x)$. 
  For each $x\in supp(c)$, let $S(x)$ be the union of all closed 2-cells  of  $X$ 
that  give rise to  $c_x$.   Since $c_x$ is a simple loop, $S(x)$ is 
 homeomorphic   to the   closed unit disk in the  Euclidean   plane.

\b{Le}\label{surface}
{ Let  $c\subset \p_T X$ be a  simple  closed geodesic  with $length(c)<4\pi$.  Then for each $x\in supp(c)$,
  there is some $r>0$   with 
  $supp(c)\cap \ol{B}(x,r)=S(x)\cap \ol{B}(x,r)$.}
\end{Le}
\b{proof}
For $x\in supp(c)$, choose $r_x$ such that 
 $\ol{B}(x,r_x)$ is isometric
 to the Euclidean cone $C_{r_x}(Link(X, x))$ over $Link(X,x)$ with radius $r_x$. 
Set $K=\ol{B}(x, r_x/2)$.
Notice $f_{x,r_x}$ represents a class in $H_2(X, X-K)$.   
By Lemma \ref{simpleloop}  $\log_x(c)$ is homotopic to a simple closed geodesic 
 in $Link(X,x)$. It then follows from the definitions of 
$S(x)$ and $f_{x,r_x}$  that 
$[f_{x,r_x}]=[S(x)\cap \ol{B}(x,r_x)]\in H_2(X, X-K) $.

For any $y\in \ol{B}(x,r_x/4)$,   by using homotopy along geodesic segments
 we see  $f_{y, r_x}$ and $f_{x, r_x}$  as maps 
  from  $(C_{r_x}(c), \p C_{r_x}(c))$   
   into the pair
 $(X,X-K)$ are homotopic. As a result, $[f_{y, r_x}]=[f_{x, r_x}]\in H_2(X,X-K)$.
 Combining with  the observation from last paragraph we see 
 $[f_{y, r_x}]=[S(x)\cap \ol{B}(x,r_x)]\in H_2(X, X-K) $. 
 It follows that for any $y\in \ol{B}(x,r_x/4)$,
   $f_{y, r_x}$ represents a  nontrivial class in $H_2(X, X-\{y\})$ if and only if 
  $S(x)\cap \ol{B}(x,r_x)$ represents a nontrivial class in $H_2(X, X-\{y\})$.
Recall $S(x)\cap \ol{B}(x,r_x)$  is homeomorphic to the closed  unit  disk
 in the  Euclidean  plane.  Now it follows easily from excision that for any $y\in \ol{B}(x,r_x/4)$, 
 $ f_{y, r_x}$ represents a nontrivial class in $H_2(X, X-\{y\})$ if and only if 
  $y\in \ol{B}(x,r_x/4)\cap (S(x)\cap \ol{B}(x,r_x))=S(x)\cap \ol{B}(x,r_x/4)$.

\end{proof}

A    similar   argument   shows the complement of $supp(c)$ in $X$ is open:

\b{Le}\label{closed}
{Let  $c\subset \p_T X$ be  a  simple  closed geodesic in 
$\p_T X$ with length strictly less than $4\pi$.   Then $supp(c)$ is a closed subset of $X$.}
\end{Le}

\b{Prop}\label{shortloop}
{Let $X$ be  a  locally compact $CAT(0)$    square complex and $c\subset \p_T X$ a  
    simple  closed geodesic in 
$\p_T X$ with length at most $2.5\pi$. Then  $supp(c)$ is a  closed 
 convex subset of $X$  and  is homeomorphic to  the plane.}
\end{Prop}
\b{proof}
Lemma \ref{surface}  implies $supp(c)$ is a 2-dimensional manifold.  
We claim $supp(c)$ is locally convex in $X$. 
  Fix  $x\in supp(c)$ and let $c_x$ be the simple loop in $Link(X,x)$ homotopic to
  $\log_x(c)$.  Let $r>0$ be the number from Lemma \ref{surface}.
  In particular  $\overline B(x,r)$  is isometric to 
  the Euclidean cone $C_r(Link(X,x))$ over the link $Link(X,x)$ with radius $r$. 
We write a point $P\in C_r(Link(X,x))$ as $P=t\xi$ where $0\le t\le r$ and $\xi\in Link(X,x)$.
 Lemma \ref{surface}
implies  $supp(c)\cap \overline B(x,r)\subset \ol{B}(x,r)$  is isometric to 
  $\{t\xi: 0\le t\le r, \xi\in c_x\}\subset C_r(Link(X,x))$.  We need to show 
$\{t\xi: 0\le t\le r, \xi\in c_x\}$ is convex in $C_r(Link(X,x))$.

 Let $A\subset Link(X,x)$ be a subset. Then it   is  not hard  to check that
 the  subset $\{t\xi: 0\le t\le r, \xi\in A\}$ is convex in $C_r(Link(X,x))$ if and only if
  for any   $a,b \in A$ with $d_x(a,b)<\pi$, $A$ contains the minimal geodesic
 $ab$ in $Link(X,x)$ from $a$ to $b$.   Since $X$ is a    $CAT(0)$  square complex 
     and the length of $c$ is at most $2.5\pi$,    $c_x$ has 
length $2\pi$ or $2.5\pi$.  Now it   is   easy to check that for any  
 $a,b \in c_x$ with $d_x(a,b)<\pi$, $c_x$ contains the minimal geodesic
 $ab$ in $Link(X,x)$ from $a$ to $b$. It follows that $supp(c)$ is locally convex in $X$.
Since $supp(c)$ is also closed, it is convex in $X$. Therefore with the induced  metric,
  $supp(c)$ is a $CAT(0)$    space. In particular $supp(c)$ is contractible. 
  As a contractible surface $supp(c)$ is homeomorphic to the plane.

\end{proof}

Let  $X$   and $c$ be  as in Proposition \ref{shortloop}.  Then  $supp(c)$  is a 
  subcomplex of $X$.  If we equip  $supp(c)$ with the path metric,  then
  $supp(c)$  is a  $CAT(0)$ square complex  and the inclusion  $supp(c)\subset  X$   
  is  an isometric embedding.   It follows that $\p_\i supp(c)$   naturally embeds 
  into $\p_\i X$.  We  identify
 $\p_\i supp(c)$ with its embedding in $\p_\i X$.


\b{Le}\label{limitset}
{Let $X$ be  a  locally compact $CAT(0)$    square complex and $c\subset \p_T X$ a  
    simple  closed geodesic in 
$\p_T X$ with length at most $2.5\pi$.
 Then $\p_\i supp(c)=c$.}
\end{Le}
\b{proof}
We first show $\p_\i supp(c)\subset c$.  
Fix a point $p\in supp(c)$   and let $\xi\in \p_\i supp(c)$. 
We need to show   $\xi\in c$.    Since  $c$ is   a  closed   subset  of  $\p_\i X$ 
    it  suffices to  show that  for any  $t\ge 0$
 there is some $\eta\in c$  with
 $\gamma_{p\xi}(t)\in p\eta$.  

Suppose there is some $t_0$ such that $\gamma_{p\xi}(t_0)\notin p\eta$  for any 
 $\eta\in c$.  Let $q=\gamma_{p\xi}(t_0)$ and  $K=\ol{B}(p, 2t_0)$. 
The proof of Lemma \ref{surface}
shows that   $f_{p,4t_0}$ and  $f_{q,4t_0}$ represent the same class in 
 $H_2(X, X-K)$.  Let $i_*: H_2(X, X-K)\ra H_2(X,X-\{q\})$ be the homomorphism induced 
 by the inclusion.  The assumption  $q\notin p\eta$  for any 
 $\eta\in c$ implies $q$ does not lie in the image of $f_{p,4t_0}$ and thus
$i_*([f_{p,4t_0}])=0\in H_2(X,X-\{q\})$.   It follows that 
$i_*([f_{q,4t_0}])=0\in H_2(X,X-\{q\})$, which implies $q\notin supp(c)$, a contradiction.

$supp(c)$ is  a  $CAT(0)$ square complex and is 
homeomorphic to the plane. 
It follows from \cite{KO} that
$\p_\i supp(c)$ is homeomorphic to a circle. Now it   is   clear that $\p_\i supp(c)=c$
  since  $\p_\i supp(c)$ is a circle embedded in 
 the   circle $c$. 
 
\end{proof}

\subsubsection{Fake Planes} \label{subfake}

The standard quarter plane   
is the first quadrant $\{(x,y)\in \E^2: x\ge 0,  \; y\ge 0\}$ with the obvious square complex structure.
A quarter plane in a   $CAT(0)$ square complex  $X$ is the image of a cellular embedding from the 
standard quarter   plane into $X$.  
 The image of the origin  is called the cone point of the quarter plane. 
The Euclidean plane  clearly  is the union of four quarter planes
  with  a  common  cone point.

\b{Def}\label{deffake}
{Let $X$ be a $CAT(0)$ square complex.  
A  subcomplex of $X$ is  a \e{fake plane}   if  it is homeomorphic to the Euclidean plane and 
is  the union of five quarter planes
 with the same cone point  and  disjoint   interiors.}
\end{Def}

It follows from the proof of Proposition \ref{shortloop}
that a fake plane in a $CAT(0)$ square complex is convex. 
 It   is   also clear if $F'$ is a  fake plane then $\p_T F'\subset \p_T X$
 is a simple   closed  geodesic  with length $2.5\pi$.

\b{Prop}\label{profake}
{Let $X $  be a $CAT(0)$ square complex.  Then there exists a simple   closed  geodesic 
 with length $\frac{5\pi}{2}$ in $\p_T X$ if and only if there  exists a fake plane in $X$.}
\end{Prop}

\b{proof} One direction is clear.
Let $c\subset \p_T X$ be a simple   closed  geodesic 
 with length $\frac{5\pi}{2}$  and $supp(c)\subset X$ its support set. 
 We shall prove that $supp(c)$ is a fake plane.

By Proposition   \ref{shortloop}  and Lemma \ref{limitset}
 $supp(c)$  is a closed convex subset of $X$  with $\p_\i supp(c)=c$.
It follows that the embedding $\p_\i supp(c)\subset \p_\i X$ is 
 an isometric embedding in the angular metric.
Since $c$ is a geodesic in $\p_T X$,   $\p_T supp(c)$ is also a geodesic 
 in the Tits metric. Therefore $\p_T supp(c)$   has length $2.5\pi$.

Now $supp(c)$ is a piecewise Euclidean $CAT(0)$ 2-complex which is homeomorphic to 
 the plane.  It follows from \cite{KO} that the difference 
$2\pi-length(\p_T supp(c))$ equals the total  curvature  on $supp(c)$.
  Since $supp(c)$ is a square complex, 2-cells are flat, edges are geodesics and only
 vertices  can contribute to curvature.   Notice  each vertex  contributes  either  0 or 
 at most $-\pi/2$. Since  $2\pi-length(\p_T supp(c))=-\pi/2$,  
we see there is exactly one vertex  $v\in supp(c)$ such that there are 5 squares 
in $supp(c)$ incident to $v$.  At every other vertex  $v'$ of $supp(c)$ there are exactly
 4 squares incident to  $v'$.  Now it   is   easy to see that  $supp(c)$ is a  fake plane.

\end{proof}

If $F'$ is a fake plane in a $CAT(0)$ square complex $X$ and $x\in F'$    is  the common cone point
 of its five quarter planes,  then  $Link(F', x)\subset Link(X,x)$ is a simple loop   consisting  of 
  5 edges. The following corollary provides a  practical  sufficient
 condition for the nonexistence of simple   closed   geodesic  with length $2.5\pi$ in 
 $\p_T X$,  especially  when $X$ admits a    cocompact action.

\b{Cor}\label{link}
{Let $X $  be a $CAT(0)$ square complex.  If no vertex link of $X$ contains a simple loop
    consisting   of   5 edges, then there is no 
 simple   closed   geodesic 
 with length $2.5\pi$ in $\p_T X$.} 
\end{Cor}

\subsection{Groups Containing Generic Elements}\label{gelement}

 Let $X $  be a $CAT(0)$ square complex,  $G$ a group acting properly 
 and cocompactly  by  cellular isometries on $X$ and $H<G$ a subgroup of $G$. 
Our goal in Sections \ref{gelement}   and \ref{secsingular}
  is to prove Theorem \ref{main}.
 In this section we consider the case 
 when the group $H$ contains  generic elements.  The main result of this section is 
Theorem \ref{ggeneric}.


\subsubsection{Parallel Set and Generic Elements}  \label{subgeneric}

Let $g$ be a hyperbolic  isometry of   $X$, and 
 $P_g=Y\times R$ its parallel set.  Since $X$ is a 2-complex,   $Y$ is a tree.

\b{Le}\label{parallel}
{Let $g\in G$ be a hyperbolic  isometry of   $X$, and 
 $P_g=Y\times R$ its parallel set.  Then the number of ends of the tree 
  $Y$  is one of the following: 0, 2, $\i$.}
\end{Le}
\b{proof} Notice $g(P_g)=P_g$ and $g_{|P_g}=(g_1,T):Y\times R  \ra Y\times R$ where
 $g_1$ is an isometry of $Y$ and $T$ is a  translation of $R$. 
Suppose  the number of ends of $Y$ is finite. Then there is an integer $n$ such that
 $g_1^n$   fixes all the ends of $Y$.  It follows that $Min(g^n)=Y'\times R$, 
 where $Y'\subset Y$ contains the  core of $Y$ (the core of $Y$ is 
the union of all geodesics between the ends of $Y$). $Y'$ has the same number of ends as $Y$.
By Theorem \ref{mincentral}, $C_{g^n}(G)$ leaves $Min(g^n)$ invariant and acts on 
$Min(g^n)$ properly and cocompactly.  It follows that $C_{g^n}(G)/<g^n>$
acts properly and cocompactly  on $Y'$. Thus the number of ends of $Y'$ equals the 
number of ends of the group $C_{g^n}(G)/<g^n>$, which takes the values 0, 1, 2 or $\i$.
It can not be 1 since any one-ended tree does not admit any proper and cocompact group
 action.
\end{proof}

We notice $g$ is of rank one if and only if $Y$ has 0 ends. 

\b{Th}\label{th6.2}
{Let $X $  be a $CAT(0)$ square complex,  $G$ a group acting properly 
 and cocompactly  by  cellular isometries on $X$, 
and 
   $g_1,  g_2\in G$  two hyperbolic isometries. 
    If   $d_T(\xi, \eta)\ge \pi$ for 
 any $\xi\in \{g_1(+\i), g_1(-\i)\}$    and  any
   $\eta\in \{g_2(+\i), g_2(-\i)\}$,  
   then  the subgroup   generated by  $g_1$  and  $g_2$   contains  a free group 
 of rank two.}
\end{Th}

\b{proof}
If one of $g_1$, $g_2$ is  of  rank one,  then the theorem follows from 
  Corollary \ref{1possibility}.  Assume    none   of  $g_1$, $g_2$    is  of rank one. 
 Then by Lemma \ref{parallel} the four points  $g_1(+\i)$,   $g_1(-\i)$,
      $g_2(+\i)$,    $g_2(-\i)$ 
  are not terminal points and the theorem follows from
    Theorems \ref{criterion}
   and \ref{pivisi}.

\end{proof}

Recall a \e{flat} in a $CAT(0)$ $2$-complex is the image of an isometric  embedding from 
 the Euclidean plane into  the $CAT(0)$ $2$-complex. 

\b{Def}\label{generic}
{Let $g$ be a hyperbolic  isometry of   $X$ and $P_g=Y\times R$ its parallel set.
We call $g$ a \e{generic element}
 if  $P_g$ is a flat and  the axes of $g$ are not parallel to any  edges in $P_g$.
 We say $g$ is a \e{special element} if $g$ is neither  generic nor of rank one.}
\end{Def}

  By  Lemma \ref{parallel}    we see 
 if $g$ is a  special element, then $g(+\i)$ and $g(-\i)$ are singular points
  in $\p_T X$.

 By Proposition \ref{specialig}   there are exactly four singular points
 on each unit circle in $\p_T X$. These four singular points are evenly spaced on the circle. 
 For a generic element $g$,        
    set 
$\alpha_g=\min\{d_T(g(+\i), s): s\in \p_T {P_g}\;\; \text{is a singular point}\}$.
 Clearly $0<\alpha_g\le \pi/4$. 
Now we are ready to state the main theorem of  Section  \ref{gelement}.

\b{Th}\label{ggeneric}
{With the assumptions   of   Theorem   \ref{main}. 
If  a subgroup $H$ of $G$  contains a generic element, then $H$ either is virtually free abelian or 
 contains a free group of rank two.}
\end{Th}

\b{remark}\label{noneed}
{Theorem \ref{ggeneric} holds without assuming there is no fake plane in $X$. 
The proof in the general case is similar but longer.}
\end{remark}

The proof of Theorem \ref{ggeneric} shall be completed in the next two sections. 
 We first need to look at the relative position of the parallel  sets of two generic
 elements. 

\b{Le}\label{vfreeab}
{Let  $h\in H$ be a generic element.  
If    $k(\p_T {P_h})=\p_T {P_h}$ 
for all $k\in H$,   then  the group $H$ is virtually free abelian.}
\end{Le}
\b{proof}
By a theorem of Leeb (\cite{L}),  $F\ra \p_T F$  
 is a  1-1  correspondence 
 between the set of flats in $X$ and the set of unit circles in $\p_T X$. 
Thus $k(\p_T {P_h})=\p_T {P_h}$   implies   $k(P_h)=P_h$.   It follows that the flat
 $P_h$ is invariant under the action of $H$ and the lemma follows from 
 Bieberbach's theorem.

\end{proof}

\b{Le}\label{3possibility}
{Let  $h, k \in H$ be  generic elements.  Then    one of the following holds:\newline
\e{(i)}   $\p_T {P_h}=\p_T {P_k}$;   in this case  $<h, k>$ is virtually free abelian;\newline
\e{(ii)}    $\p_T {P_h}\cap \p_T {P_k}=\phi$;\newline
\e{(iii)}   $\p_T {P_h}\not=\p_T {P_k}$   and $\p_T {P_h}\cap \p_T {P_k}\not=\phi$; in this case  
  $<h^2, k^2>$   is   free  of rank two.}
\end{Le}
\b{proof}
If $\p_T{P_h}=\p_T{P_k}$, then  $P_h=P_k$  and   both $h$ and $k$ leave  invariant 
    the flat $P_h=P_k$. It follows from 
 Bieberbach's theorem    that 
$<h, k>$ is virtually free abelian.
Now suppose $\p_T{P_h}\not=\p_T{P_k}$ and $\p_T {P_h}\cap \p_T {P_k}\not=\phi$. 
Then for some $\epsilon>0$, the intersection $N_\epsilon(P_h)\cap N_\epsilon(P_k)$
 is nonempty and unbounded.  
Notice  $Min(h^2)=P_h$  and $Min(k^2)=P_k$.  By Theorem \ref{mincentral}, $C_{h^2}(G)$
   leaves  invariant  $P_h$ and acts on $P_h$   properly  and cocompactly.  
  Similarly  for $C_{k^2}(G)$
    and   $P_k$.   
By Theorem \ref{swenson},  the group
 $C_{h^2}(G)\cap C_{k^2}(G)$ acts     properly    and   
 cocompactly on $N_\epsilon(P_h)\cap N_\epsilon(P_k)$.
  Since   the  set  $N_\epsilon(P_h)\cap N_\epsilon(P_k)$  is  unbounded, 
       a theorem of E. Swenson (\cite{Sw})  implies 
  there exists 
  a hyperbolic isometry $g\in  C_{h^2}(G)\cap C_{k^2}(G)$.  
Then both $h^2$ and $k^2$ commute with $g$    and  leave  $P_g$ invariant.  
$P_g$ splits $P_g=Y\times R$ where $Y$ is a tree. 
   Since   $h^2$ and $k^2$ are generic elements   and  $\p_T{P_h}\not=\p_T{P_k}$,  the axes 
  of $h^2$ and $k^2$ are not parallel to the axes of  $g$.   It follows that 
there are geodesics 
 $c_1, c_2\subset Y$ such that $P_h=c_1\times R$ and $P_k=c_2\times R$.  
  $\p_T{P_h}\not=\p_T{P_k}$   implies  $c_1$ and $c_2$ share   at   most   one end.
If  $c_1$ and $c_2$ share exactly one end, then  $P_h\cap P_k$ is a flat half-plane. 
 On the other hand  by Theorem \ref{swenson} the group 
$C_{h^2}(G)\cap C_{k^2}(G)$ acts cocompactly on  $P_h\cap P_k$. This is a contradiction since 
 no group acts cocompactly on a flat half-plane.  Therefore $c_1$ and $c_2$ share no end. 
  In this case  it   is   easy to see $<h^2, k^2>$   is  free  of rank two  by looking at   the  
 induced actions of  $h^2$  and  $k^2$  on $Y$. 

\end{proof}

\subsubsection{When the Minimal Set  Contains Singular Points }\label{containsspe}

We  prove Theorem \ref{ggeneric} in this section and  Section \ref{containsno}.

We assume $H$  neither    is   virtually 
free abelian nor contains a free group of rank two  and will 
   derive a contradiction from this.
By Corollary \ref{1possibility} we conclude  $H$ contains no rank one isometries.
Let $\Lambda(H)\subset \p_\i X$ be the limit set of $H$ and $M\subset \Lambda(H)$ a minimal set of $H$.
Theorem \ref{dis} implies 
 $d_T(\xi,m)\le \pi$ for any $\xi\in \Lambda(H)$  and any $m\in M$.
  In particular, 
$d_T(m, g(+\i) )\le \pi$, $d_T(m, g(-\i) )\le \pi$
 for any hyperbolic isometry  $g\in H$  
and any $m\in M$ (since $g(+\i), g(-\i) \in  \Lambda(H)$).

Since $M$ is closed and 
  $H$-invariant, Theorem \ref{hyperbolicf}  implies for any hyperbolic isometry $g\in H$, 
  $M\cap \p_T {P_g}\not=\phi$.   By last paragraph  $g$  is  not 
 a rank one isometry.  Lemma \ref{parallel}
  implies     a point in  $\p_T {P_g}\supset M\cap \p_T {P_g}$ is either 
 singular or regular.

Let $h\in H$ be a generic element. By Lemma \ref{vfreeab} there is some $k\in H$ such that
$k(\p_T{P_h})\not=\p_T{P_h}$.  Let $h'=khk^{-1}$. Then $h, h'\in H$ 
are  two generic elements  with $\p_T{P_h}\not=\p_T{P_{h'}}$.   By Lemma \ref{3possibility}
$\p_T{P_h}\cap \p_T{P_{h'}}=\phi$. Also notice $\alpha_h=\alpha_{h'}$.

In this section we prove 
 Theorem   \ref{ggeneric}  when the minimal set $M$ contains  singular points:

\b{Prop}\label{genespe}
{Theorem \ref{ggeneric} holds if  there is a  generic element  $h_0\in H$ such that 
$M\cap \p_T {P_{h_0}}$
 contains  a  singular point.}

\end{Prop}


Suppose  $h\in H$ is  a generic element  and $s\in M-\p_T {P_h}$  is 
  a singular point  but not a terminal point.
Then $d_T(h(-\i),s)\le \pi$.  Since $h(-\i)$ is    a regular point and $s$ is a singular point
Proposition \ref{specialig} implies 
$d_T(h(-\i),s)< \pi$  and  there is  at most one singular point  in the interior of  the 
 geodesic segment $h(-\i)s$.  On the other hand, since $s\notin\p_T {P_h}$ 
 there is 
 some  $\xi$  in the interior of $h(-\i)s$ 
 such that  $h(-\i)s\cap \p_T {P_h}=h(-\i)\xi$.
 $\xi$ is a branch point, and by Proposition  \ref{branch}
it is also a singular point.
 Notice $d_T(\xi, s)=\pi/2$.

\b{Le}\label{5.2lemma1}
{Suppose $h\in H$ is  a generic element, $s\in M-\p_T {P_h}$ is  a singular point 
   but not a terminal point  
   and $\xi$   is   the unique singular point in the interior of $h(-\i)s$. 
  Denote by $\eta$  the unique point  on $\p_T {P_h}$ with   
$d_T(\eta, h(-\i))=d_T(s, h(-\i))$   and  $d_T(\eta, \xi)=\pi/2$.
Then $\eta\in M$   and is a singular point. 
 Similarly, if $\xi'$ is 
the unique singular point in the interior of $h(+\i)s$   and 
$\eta'$ is  the unique point  on $\p_T {P_h}$ with  
$d_T(\eta', h(+\i))=d_T(s, h(+\i))$  and  $d_T(\eta', \xi')=\pi/2$.
Then $\eta'\in M$ and  is a singular point.}
\end{Le}
\b{proof}
We prove the claim on $\eta$, the proof for $\eta'$ is similar.  
   By Proposition \ref{specialig} 
$\eta$ is a   singular point. 
 By replacing $h$ with $h^2$
  if necessary 
   we may assume $Min(h)=P_h$.    In particular $h(\xi)=\xi$. 
Let $s'\in \p_\i X$ be an accumulation point (in the cone topology) of the set 
$\{h^i(s): i\ge 1\}$. We see  $s'\in M$  since $s\in M$ and   $M$ is 
$H$-invariant and closed. 
 Notice $d_T(h(-\i), s)=\angle_T(h(-\i),s)$ since 
$d_T(h(-\i),s)< \pi$ (Proposition \ref{Titsb}). 
By Theorem \ref{hyperbolicf} we have  $s'\in\p_T {P_h}$ and   $d_T(h(-\i), s')=d_T(h(-\i), s)$.
On the other hand, $h(\xi)=\xi$ and a subsequence of $\{h^i(s):i\ge 1\}$
 converges to $s'$ in the cone topology. 
Proposition \ref{Titsb}(iii)    implies 
 $d_T(\xi, s')\le d_T(\xi,s)=\pi/2$. It now follows that $d_T(\xi, s')=\pi/2$, 
  that is $s'=\eta$.

\end{proof}

Lemma \ref{5.2lemma1}    and the assumption in  Proposition  \ref{genespe}   imply   
 for any generic element $h\in H$,
   $\p_T {P_h}\cap M$ contains a singular point.

\b{Le}\label{5.2lemma2}
{Let $h\in H$ be a generic element, and $s\in M-\p_T {P_h}$ a singular
    but  not a terminal point.
 Let  $\xi_1$,   $\xi_2$ be the unique singular points in the interior of 
 $sh(+\i)$ and $sh(-\i)$ respectively. 
  Then $sh(+\i)\cap sh(-\i)=\{s\}$  and  $d_T(\xi_1, \xi_2)=\pi$.}
\end{Le}
\b{proof}   
 Suppose  $sh(+\i)\cap sh(-\i)=s\xi$ is a nontrivial segment.  
Then $\xi$ is a singular point since it is a  branch point. 
 Proposition \ref{specialig}   implies  $d_T(s, \xi)=\pi/2$.
Since 
$d_T(s, h(+\i))<\pi$   and   $d_T(s, h(-\i))<\pi$,    we   have 
$d_T(\xi, h(+\i))<\pi/2$  and   $d_T(\xi, h(-\i))<\pi/2$. It follows that
 $d_T(h(+\i), h(-\i))<\pi$, a contradiction. 

By the paragraph preceding Lemma \ref{5.2lemma1}, 
   $\xi_1,  \xi_2\in \p_T {P_h}$.  If $d_T(\xi_1, \xi_2)=\pi/2$,
  then $\xi_1s\cup s\xi_2\cup\xi_2\xi_1$ is a simple  closed  geodesic  with length $1.5\pi$
    in the 
 $CAT(1)$ space $\p_T X$, a contradiction.

\end{proof}

 Now we can complete the proof of Proposition \ref{genespe}.

\noindent
{\bf{Proof of Proposition \ref{genespe}.}}
There are two generic elements $h, k\in H$ such that 
$\p_T{P_h}\cap \p_T{P_{k}}=\phi$   and  $\alpha_h=\alpha_{k}$. 
   After replacing 
 $h$ and $ k$  with $h^2$ and $k^2$ if necessary
  we may assume $Min(h)=P_h$ and $Min(k)=P_k$. 
By   Theorem \ref{hyperbolicf}
  the fixed point sets of $h$ and $k$ in  $\p_\i X$ are 
  $\p_T{P_h}$ and $\p_T P_{k}$ respectively.

Let  $\eta_3\in \p_T {P_k}\cap M$  be a  singular point, 
 and $\xi_1$, $\xi_2$ be singular points in the interior of 
$h(-\i)\eta_3$ and $h(+\i)\eta_3$  respectively. Then $\xi_1$, $\xi_2\in \p_T {P_h}$.
 Lemma \ref{5.2lemma2}  implies 
  $d_T(\xi_1, \xi_2)=\pi$.      By   Lemma \ref{5.2lemma1}  
 the other two singular points (denoted  by  $\eta_1$, $\eta_2$) on $\p_T {P_h}$ lie
 in $M$.  Let $\xi_3$, $\xi_4$ be the  two singular points  on $\p_T {P_k}$  
 that have distance $\pi/2$ from $\eta_3$, and 
$\eta_4$  the singular point on $\p_T {P_k}$  with $d_T(\eta_3, \eta_4)=\pi$.
   Since by assumption there is no fake plane in $X$,  Proposition \ref{profake}
   implies  there is no simple closed geodesic with length $2.5\pi$ in $\p_T X$. 
    Of course  there is also  no simple   closed geodesic  with
 length $1.5\pi$ in  the $CAT(1)$ space   $\p_T X$. 
  By considering the unique singular points  in the interior of the geodesic segments 
  from $k(+\i)$ and $k(-\i)$  to $\eta_i$ ($i=1,2$),  we see 
$d_T(\eta_i, \xi_j)=\pi/2$ for $i=1,2$ and $j=3,4$.
 Now Lemma \ref{5.2lemma1} applied to the generic 
    element $k$ or $k^{-1}$ and singular point  $\eta_1\in M-\p_T {P_k}$ implies 
  $\eta_4\in M$.    Again by considering the unique 
singular points  in the interior of  $h(+\i)\eta_4$ and $h(-\i)\eta_4$,  we conclude 
 $d_T(\xi_i,\eta_j)=\pi/2$ for any $i,j\in\{1,2,3,4\}$.

We may assume 
            $d_T(k(+\i), \eta_3)<\pi/2$
  by replacing  $k$ with $k^{-1}$   if necessary.
    Note the choice  of $\xi_2$ implies $d_T(h(+\i), \xi_2)<\pi/2$.
  If   $d_T(h(+\i), \xi_2)+d_T(\eta_3, k(+\i))=\pi/2$, then it   is   easy to check that 
  $d_T(\xi, \eta)=\pi$ for any $\xi\in \{h(+\i),h(-\i)\}$ and any 
$\eta\in \{k(+\i),  k(-\i)\}$.  Now Theorem   \ref{th6.2}      implies that 
 $<h, k>$ contains a free group of rank two.

 If   $d_T(h(+\i), \xi_2)+d_T(\eta_3, k(+\i))\not=\pi/2$,  we consider 
$hkh^{-1}$ and $k$ instead of $h$ and $k$.  
 Notice $\p_T {P_{hkh^{-1}}}=h(\p_T {P_k})\not=\p_T {P_k}$.
       If $\p_T {P_{hkh^{-1}}}\cap \p_T {P_k}\not=\phi$, then by Lemma \ref{3possibility}
      $<hkh^{-1}, k>$    contains   a free group of rank two and we are done. Assume 
  $\p_T {P_{hkh^{-1}}}\cap \p_T {P_k}=\phi$. 
Since $d_T(\eta_1, \xi_3)=d_T(\eta_1, h(\xi_3))=\pi/2$,  we see $d_T(\xi_3, h(\xi_3))=\pi$.
  The argument on $h$ and $k$ shows  that $d_T(\xi_i, h(\eta_j))=d_T(h(\xi_i), \eta_j)=\pi/2$
 for any $i,j\in \{3,4\}$.   Notice now we have 
 $d_T(hkh^{-1}(+\i), h(\xi_4))+d_T(\eta_3, k(+\i))=\pi/2$. 
Now it   is  easy to see $d_T(\xi, \eta)=\pi$ for any $\xi\in \{hkh^{-1}(+\i),hkh^{-1}(-\i)\}$ and any 
$\eta\in \{k(+\i),  k(-\i)\}$   and Theorem   \ref{th6.2}  completes the proof.

\qed

\subsubsection{When the Minimal Set  Contains   Regular Points }\label{containsno}

In this section we prove 
Theorem   \ref{ggeneric}  when the minimal set $M$ contains  regular points:


\b{Prop}\label{5.3pro}
{Theorem \ref{ggeneric} holds if  there is a  generic element  $h_0\in H$ such that 
$M\cap \p_T {P_{h_0}}$
 contains  a  regular point.}

\end{Prop}

We notice  Theorem \ref{ggeneric} follows from Propositions \ref{genespe} and \ref{5.3pro}.

The proof of Proposition \ref{5.3pro} is similar to that of Proposition \ref{genespe}.
We  assume $H$    neither   is  virtually 
free abelian nor contains a free group of rank two  and will 
   derive a contradiction from this. 
  As in Section \ref{containsspe} we see:
there are two generic elements $h, k\in H$ such that 
$\p_T{P_h}\cap \p_T{P_{k}}=\phi$   and  $\alpha_h=\alpha_{k}$;
$d_T(m, g(+\i) )\le \pi$   and  $d_T(m, g(-\i) )\le \pi$ 
 for any hyperbolic isometry  $g\in H$  
and any $m\in M$.

\b{Le}\label{twopairs}
{Let $h\in H$ be a generic element and $s_1, s_2\notin \p_T {P_h}$ singular 
           but  not  terminal  points with
 $d_T(s_1,s_2)=\pi/2$. Let $\xi_0\in \{h(+\i), h(-\i)\}$  
   and     $s_3, s_4$ the two singular points 
 on  $\p_T {P_h}$ that have distance less than $\pi/2$ from $\xi_0$.
If $m\in M$ lies in the interior of $s_1s_2$,   then there are 
 $\xi\in \{s_3, s_4\}$ and $\eta\in \{s_1, s_2\}$ such that $d_T(\xi,\eta)=\pi/2$.}

\end{Le}
\b{proof}  We have   $d_T(\xi_0, m)\le \pi$ 
    by the paragraph   preceding the lemma.
  Let $\sigma$ be a minimal geodesic from $m$ to $\xi_0$.
 Since any branch point is a singular point, there are 
 $\xi\in \{s_3, s_4\}$ and $\eta\in \{s_1, s_2\}$ such that  $\xi$ and 
$\eta$ lie in the interior of $\sigma$. Clearly $\xi\not=\eta$. 
 On the other hand, 
$d_T(m, \xi_0)\le \pi$   and  Proposition \ref{specialig}
    imply   there are at most two singular points 
 in the interior of  $\sigma$. Therefore $d_T(\xi,\eta)=\pi/2$.

\end{proof}

\b{Le}\label{inside}
{Let $h\in H$ be an arbitrary  generic element, $\xi=h(+\i)$ or $h(-\i)$, and 
 $s_1$, $s_2$ the two singular points on $\p_T{P_h}$ with $d_T(\xi, s_i)<\pi/2$.
  Then the interior of $s_1s_2$ contains a point of the minimal set $M$.}
\end{Le}  

\b{proof} 
We may assume $Min(h)=P_h$ by replacing $h$ with $h^2$ if necessary.
Since we assume $H$  neither    is  virtually 
free abelian nor contains a free group of rank two, 
   Lemmas \ref{vfreeab} and \ref{3possibility} imply that 
there is a generic element  $k\in H$ with  $\p_T {P_h}\cap \p_T {P_k}=\phi$.

Notice  $M\cap \p_T {P_k}\not=\phi$.  By Proposition  \ref{genespe}  we see 
 $M\cap \p_T {P_k}$ contains no singular points.
    Let $m\in M\cap \p_T {P_k}$   be   a regular point, and 
  $s_1'$, $s_2'$   the two singular points on  $\p_T {P_k}$
 with $d_T(m, s_i')<\pi/2$. Also let $s_3$, $s_4$ be the two 
singular points on  $\p_T {P_h}$  other than $s_1$    and  $s_2$, and
  $\xi'\in \{ h(+\i), h(-\i)\}$  but $\xi'\not=\xi$.    
Notice $d_T(\xi', m)\le \pi$. 
Let $\sigma$ be a minimal geodesic from $m$ to $\xi'$.
 Lemma \ref{twopairs} implies there are 
 $\eta_1\in \{s_1', s_2'\}$, $\eta_2\in\{s_3, s_4\}$ such 
that $\eta_1$ and $\eta_2$ are the only singular points 
in the interior of $\sigma$. Notice $\pi/2<d_T(\eta_2, m)<\pi$
  and $d_T(\xi', m)=d_T(\xi', \eta_2)+d_T(\eta_2, m)$.

   If   $\xi=h(+\i)$  then we define $m'$   to  be an accumulation 
point of  $\{h^i(m): i\ge 1\}$, otherwise,   $\xi=h(-\i)$  and we 
define $m'$   to  be an accumulation 
point of   the set $\{h^{-i}(m): i\ge 1\}$.   In either case, $m'\in M\cap \p_T {P_h}$.
   It  follows from     Theorem \ref{hyperbolicf}   that  $d_T(m', \xi')=d_T(m, \xi')$.
  The proof of    Lemma \ref{5.2lemma1}
   shows  
$$d_T(\xi', m')=d_T(\xi', \eta_2)+d_T(\eta_2, m').$$
    Consequently   $\eta_2$   
   lies on a minimal geodesic from $\xi'$ to $m'$  and 
 $d_T(\eta_2, m')=d_T(\eta_2, m)$.
  Since $\pi/2<d_T(\eta_2, m)<\pi$, $m'$ lies in the interior of $s_1s_2$.

\end{proof}

   Since we assume there is no fake plane in $X$  
Proposition \ref{profake}   implies   there is no simple   closed  geodesic  with length 
 $2.5\pi$ in $\p_T X$. 

\b{Le}\label{piok}
{Let $h,k\in H$ be  generic elements  with   $\p_T {P_h}\cap \p_T {P_k}=\phi$
  and $\alpha_h=\alpha_k$.  If there are $\xi_0\in \{h(+\i), h(-\i)\}$,
   $\eta_0\in \{k(+\i), k(-\i)\}$ with $d_T(\xi_0, \eta_0)=\pi$,  then 
  $d_T(\xi, \eta)=\pi$ for any $\xi\in \{h(+\i), h(-\i)\}$,
   $\eta\in \{k(+\i), k(-\i)\}$.  In particular,  $<h,k>$ contains a 
free group of rank two.}
\end{Le}
\b{proof}
Let $\sigma$ be a minimal geodesic from $\xi_0$ to $\eta_0$,   and
 $\xi_1\in \p_T {P_h}$, $\eta_1\in \p_T {P_k}$ the only two singular points on $\sigma$. 
 Clearly $d_T(\xi_1, \eta_1)=\pi/2$   
 and  $d_T(\xi_0,\xi_1)+d_T(\eta_1,\eta_0)=\pi/2$.   Let 
 $\xi_0'$ and $\eta_0'$  be   defined by  $\{\xi_0, \xi_0'\}=\{h(+\i), h(-\i)\}$,
$\{\eta_0, \eta_0'\}=\{k(+\i), k(-\i)\}$.  We denote the  singular points
on  $\p_T {P_h}$ by  $\xi_1, \eta_3, \xi_2, \eta_4$ such that
 the six points  $\xi_1, \xi_0, \eta_3, \xi_2, \xi_0',  \eta_4$  are in cyclic order on
$\p_T {P_h}$.  Similarly we denote the singular points on 
$\p_T {P_k}$ by  $\eta_1, \xi_3, \eta_2, \xi_4$ such that
 the six points    $\eta_1,\eta_0,  \xi_3, \eta_2, \eta_0', \xi_4$ are in cyclic order on
$\p_T {P_k}$.

 Lemma \ref{inside}   implies   each of the four segments $\xi_1\eta_3$, 
$\xi_2\eta_4$, $\eta_1\xi_3$, $\eta_2\xi_4$  contains    a 
point of the minimal set $M$ in its interior.
 Lemma \ref{twopairs} applied to the two segments  $\xi_1\eta_3$ and $\eta_2\xi_4$ shows
    that
 at least one of the following holds: (1) $d_T(\xi_1, \eta_2)=\pi/2$; \linebreak
(2) $d_T(\xi_1, \xi_4)=\pi/2$;
  (3) $d_T(\eta_3, \eta_2)=\pi/2$;
 (4) $d_T(\eta_3,\xi_4)=\pi/2$.    (2) can  not hold since otherwise 
$\xi_1\xi_4\cup \xi_4\eta_1\cup\eta_1\xi_1$ is a simple   closed  geodesic  with length $1.5\pi$.
(3) can   not hold since otherwise  
$\eta_1\xi_1\cup \xi_1\eta_3\cup \eta_3\eta_2\cup \eta_2\xi_4\cup \xi_4\eta_1$ 
is a simple   closed geodesic  with length   $2.5\pi$.  So either (1) or (4) holds.


Suppose (1)  holds.  Then $d_T(\xi_0,\eta_0')=\pi$  since 
  $d_T(\xi_0, \xi_1)+d_T(\xi_1, \eta_2)+d_T(\eta_2, \eta_0')
=d_T(\xi_0, \xi_1)+\pi/2+d_T(\eta_1, \eta_0) =\pi$.   Apply   Lemma  \ref{twopairs} to
 the two segments  $\eta_4\xi_2$  and  $\eta_1\xi_3$.  
The argument in the preceding paragraph shows   either
  $d_T(\eta_4, \xi_3)=\pi/2$  or  $d_T(\xi_2,\eta_1)=\pi/2$.  
  If $d_T(\eta_4, \xi_3)=\pi/2$,  then  
$d_T(\xi_0', \eta_4)+d_T(\eta_4, \xi_3)+d_T(\xi_3, \eta_0)
=d_T(\xi_0', \eta_4)+\pi/2+d_T(\xi_1, \xi_0)=\pi$. Thus $d_T(\xi_0', \eta_0)=\pi$.
 Similarly we show  $d_T(\xi_0', \eta_0)=\pi$ if $d_T(\xi_2,\eta_1)=\pi/2$.
  By applying Lemma  \ref{twopairs} to the two segments
$\xi_2\eta_4$ and $\xi_4\eta_2$  and using a similar argument we   also  
conclude $d_T(\xi_0', \eta_0')=\pi$.

The proof when (4) holds is similar.

\end{proof}

\b{Le}\label{dlesspi}
{Let $h,k\in H$ be  generic elements  with   $\p_T {P_h}\cap \p_T {P_k}=\phi$
  and $\alpha_h=\alpha_k$.  
  Then  
  $<h,k>$ contains a 
free group of rank two.}
\end{Le}
\b{proof}   We may assume $Min(h)=P_h$  after replacing $h$ with $h^2$ if necessary.
Set $\eta_0=k(+\i)$. 
By    Lemmas \ref{twopairs}   and  \ref{inside}  we  see there are   singular
  points $\xi_1\in \p_T {P_h}$,  $\eta_1, \eta_2\in \p_T {P_k}$ with
 $d_T(\xi_1, \eta_1)=\pi/2$   and  $d_T(\eta_0, \eta_i)<\pi/2$ ($i=1,2$).  

Consider $hkh^{-1}$, $k$ and     $h(\p_T {P_k})$,    $\p_T {P_k}$.  Notice 
 $h(\p_T {P_k})=\p_T {P_{hkh^{-1}}}$.     It follows from   
  $\p_T {P_h}\cap \p_T {P_k}=\phi$   that   $\p_T {P_k}\not= \p_T {P_{hkh^{-1}}}$.
By  Lemma   \ref{3possibility}
  we may assume $\p_T {P_k}\cap \p_T {P_{hkh^{-1}}}=\phi$.
Apply Lemma \ref{twopairs} to the two segments
$\eta_1\eta_2$ and   $h(\eta_1)h(\eta_2)$.  Since $d_T(\xi_1, \eta_1)=d_T(\xi_1, h(\eta_1))=\pi/2$,
 by using the fact that there is no simple   closed geodesic  with length $1.5\pi$ or  $2.5\pi$ 
  in $\p_T X$  we
 see either $d_T(h(\eta_2), \eta_1)=\pi/2$ or  $d_T(h(\eta_1), \eta_2)=\pi/2$.
  In both cases we  conclude that $d_T(\eta_0, h(\eta_0))=\pi$.    Recall $\eta_0=k(+\i)$
  and notice  $h(\eta_0)=hkh^{-1}(+\i)$.  
   Now Lemma  
\ref{piok} applied to $hkh^{-1}$ and $k$ implies that $<hkh^{-1},k>$ contains a 
free group of rank two.

\end{proof}

 Proposition \ref{5.3pro} now follows from Lemma \ref{dlesspi}.

\subsection{Groups Containing Special Elements} \label{secsingular}

In Section \ref{gelement} we studied subgroups containing generic elements.  In this section we 
 consider subgroups containing special elements.

\b{Th}\label{gsingular}
{With the assumptions   of   Theorem   \ref{main}. 
If  a subgroup $H$ of $G$  contains a special  element, then $H$ either is virtually free abelian or 
 contains a free group of rank two.}

\end{Th}

We notice Theorem \ref{main} follows  from   Corollary \ref{1possibility}   and 
Theorems \ref{ggeneric}   and \ref{gsingular}.

 As in Section \ref{gelement} we will prove Theorem \ref{gsingular} by assuming 
$H$   neither    is   virtually 
free abelian nor contains a free group of rank two and then deriving a contradiction from this. 
Let $M$ be a minimal set of $H$.  Then $d_T(\xi, m)\le \pi$ for
 any   $\xi\in \Lambda(H)$ and any $m\in M$.
We again consider two cases depending on whether
$M$ contains any singular points.

\subsubsection{When the Minimal Set Contains  Regular Points} \label{subsecnos}

In this section we prove Theorem \ref{gsingular}    when  
$M$ contains    regular points:

\b{Prop}\label{nonspesin}
{Theorem \ref{gsingular}  holds if    there   exists a  special  element  $h_0\in H$ such that 
  $M\cap \p_T {P_{h_0}}$   contains   a  regular point.}

\end{Prop}

  We recall  that if $h\in H$ is a special element,
 then  $h(+\i)$ and $h(-\i)$ are singular points. 
  Notice in the proof of the following four lemmas we 
do   not   have to assume that $M$ contains regular points. 

\b{Le}\label{6.1lemma1}
{Let $h, k\in H$ be two  special elements   such that 
the four points
 $h(+\i)$, $h(-\i)$, $k(+\i)$, $k(-\i)$ are all distinct.
 If there is some $\xi\in\p_T X$ such that
  $d_T(\xi, \eta)=\pi/2$   whenever   $\eta$ is  one  of the four points
 $h(+\i)$, $h(-\i)$, $k(+\i)$, $k(-\i)$,  then
   $<h,k>$ contains a free group of rank two.}
\end{Le}
\b{proof}
Notice $d_T(\xi', \eta')=\pi$ for any $\xi'\in \{h(+\i), h(-\i)\}$, 
 $\eta'\in \{k(+\i), k(-\i)\}$.

\end{proof}

\b{Le}\label{6.1lemma2}
{Let $h, k\in H$ be two special elements.  If 
$d_T(h(+\i), k(+\i))=\pi/2$, $d_T(h(+\i), k(-\i))=\pi/2$ and $h$ fixes at least 
one of $k(+\i)$, $k(-\i)$,
  then  $H$  either is virtually free abelian or contains 
   a free group of rank two.}  

\end{Le}

\b{proof}
The assumption implies $\p_T Min(h)\cap \p_T Min(k)\not=\phi$.  There 
 exists $\epsilon>0$ such that 
$N_\epsilon(Min(h))\cap N_\epsilon(Min(h))$ is unbounded.
Theorems  \ref{mincentral}  and \ref{swenson}   imply $C_h(G)\cap C_k(G)$   acts properly and 
  cocompactly  on  $N_\epsilon(Min(h))\cap N_\epsilon(Min(h))$.
       There is   
  a hyperbolic isometry (see proof of Lemma \ref{3possibility})
 $g\in C_h(G)\cap C_k(G)$. Both $h$ and $k$ commute with $g$. 
It follows that $h(P_g)=P_g$,  $k(P_g)=P_g$.  $h$ is not of rank one implies  $g$ is not of rank one. 
 Thus $g$ is either generic or special. 
If $g$ is generic, then $P_g$ is a flat. In this case, $<h,k>$ acts cocompactly
 on $P_g$, so contains a generic element and the lemma follows from  Theorem \ref{ggeneric}.

Suppose $g$ is a special element. $P_g$ splits $P_g=Y\times R$ where $Y$ is a tree.
The assumption implies one of $h$, $k$, 
  say $h$ has an axis of the form $\{y_0\}\times R$ ($y_0\in Y$), and
 $k$ has an axis of the form 
 $c\times \{t_0\}$ ($t_0\in R$) where $c\subset Y$ is a geodesic in $Y$. 
  Notice $h(c\times \{t_0\})$   is an axis of $hkh^{-1}$  and has the form
$h(c\times \{t_0\})=c'\times \{t_1\}$  for a geodesic $c'\subset Y$ and some $t_1\in R$.
Recall in a discrete group, the axes of any two hyperbolic isometries 
can not share exactly one endpoint.  
Thus $c\times \{t_0\}$   and  $c'\times \{t_1\}$ 
either do  not    share any endpoints or 
 have the same endpoints. If they do  not   share any endpoints, then 
 $<k, hkh^{-1}>$ is a free group of rank two.  If they have the same endpoints,
 then   $c'=c$     
 since $Y$ is a tree. 
 In this case both $h$ and $k$ leave the flat $c\times R$ invariant and    $<h,k>$   
 acts cocompactly on 
 $c\times R$ and therefore contains generic elements.

\end{proof}

\b{Le}\label{6.1lemma3}
{Let $h, k\in H$ be two special elements.  If 
$d_T(h(+\i), k(+\i))=\pi/2$, $d_T(h(+\i), k(-\i))=\pi/2$, then 
$H$  either is virtually free abelian or contains 
   a free group of rank two.} 

\end{Le}
\b{proof}
 If  $\{k(+\i), k(-\i)\}\cap \{h(k(+\i)), h(k(-\i))\}\not=\phi$, 
   then  either  there is   some   $n\in \{1,2,3\}$  so that 
$h^n$ fixes  one of  $k(+\i)$,  $k(-\i)$,  or  the two sets 
$\{k(+\i), k(-\i)\}$ and   $ \{h^2(k(+\i)), h^2(k(-\i))\}$  are disjoint.
 In the former case the    
lemma follows from Lemma  \ref{6.1lemma2}.

 Now we may assume  $\{k(+\i), k(-\i)\}\cap \{h^i(k(+\i)), h^i(k(-\i))\}=\phi$ for $i=1$ or $2$.
 It follows from Lemma \ref{6.1lemma1}  that 
 $<k, h^ikh^{-i}>$ contains a free group of rank two
  since 
 $d_T(h(+\i), \xi)=\pi/2$ 
 for any $\xi\in \{k(+\i), k(-\i), h^i(k(+\i)), h^i(k(-\i))\}$.

\end{proof}

\b{Le}\label{6.1lemma4}
{ If  
$\{h(+\i), h(-\i)\}\cap \{k(+\i), k(-\i)\}\not=\phi$
for any two special elements $h,k\in H$, then $H$ is 
   virtually  infinite cyclic.}
\end{Le}

\b{proof}
Let $k\in H$ be  a fixed  special element.  Then for any $h\in H$,
  $hkh^{-1}$ is also a special element.  The assumption implies that
   the axes of $k$ and 
 $hkh^{-1}$ are parallel.  Thus $h(P_k)=P_{hkh^{-1}}=P_k$ for any $h\in H$.
  $P_k$ splits as $Y\times R$ where $Y$ is a tree.        
After passing to an  index two subgroup if necessary, each $h\in H$ acts on
$P_k=Y\times R$  as $h=(h', T)$ where   $h'$ is an isometry of $Y$  and 
   $T$ is a translation on $R$.  Then 
 $h\ra T$ defines a homomorphism from  $H$  onto an infinite cyclic group.
  We notice the image of the homomorphism  is indeed infinite cyclic since $H$ is a group
 of cellular isometries.  
 Let $H_0$ be the kernel of this homomorphism.  There is no hyperbolic isometry
 in $H_0$ since any hyperbolic isometry in $H_0$ is special and   the assumption implies 
    its   axes should be parallel to $\{y\}\times R$.  So $H_0$ is a torsion group acting
 properly on the tree $Y$. It follows that $H_0$ is a finite group.

\end{proof}


\noindent
{\bf{Proof of Proposition \ref{nonspesin}.}}
  Since  $M\cap \p_T {P_{h_0}}$   contains   a  regular point,  
  Theorem \ref{hyperbolicf}
  and      Proposition   \ref{specialig}  imply  that  
$M\cap \p_T {P_{h}}$   contains   a  regular point  for any  
special element $h\in  H$.  
By Lemma \ref{6.1lemma4}  we may assume there are  special elements 
  $h,k\in H$  with $\{h(+\i), h(-\i)\}\cap \{k(+\i), k(-\i)\}=\phi$. 
  Let $m\in M\cap \p_T {P_k}$ be a regular point.
By replacing $k$ with $k^{-1}$ if   necessary  we may assume 
$d_T(m, k(+\i))<\pi/2$.    Let 
 $s\in \p_T {P_k}$   be 
 the   unique  singular point   in  the interior  of   $mk(-\i)$. 
 


Notice  $d_T(m, h(+\i)), d_T(m, h(-\i))<\pi$    and  
  recall a branch point in $\p_T X$ is a singular point.  
 First suppose the intersection 
 $mh(+\i)\cap mh(-\i)$ is a nontrivial segment. 
Then $mh(+\i)\cap mh(-\i)=ms$ or $mk(+\i)$.   If $mh(+\i)\cap mh(-\i)=ms$,
  then $d_T(s, \xi)=\pi/2$ for
 any $\xi\in \{h(+\i), h(-\i),  k(+\i), k(-\i)\}$.  In this case 
 the proposition   follows from Lemma \ref{6.1lemma1}.
If $mh(+\i)\cap mh(-\i)=mk(+\i)$, then 
$d_T(k(+\i), h(+\i))=d_T(k(+\i), h(-\i))=\pi/2$ and the 
proposition follows from Lemma \ref{6.1lemma3}.

Now suppose $mh(+\i)\cap mh(-\i)=\{m\}$. In this case, $h(+\i)m\cup mh(-\i)$ is a geodesic
  between two singular points and so has length an integral multiple of $\pi/2$.
Since $d_T(m, h(+\i)), d_T(m, h(-\i))<\pi$, the length of $h(+\i)m\cup mh(-\i)$
 is either $\pi$ or $1.5\pi$.  The length of $h(+\i)m\cup mh(-\i)$
can not be $1.5\pi$ since $d_T(h(+\i), h(-\i))=\pi$ and there 
is no   simple closed geodesic  with length $2.5\pi$ in  $\p_T X$.  Therefore 
the length of $h(+\i)m\cup mh(-\i)$
 is  $\pi$. It follows that   either $s\in \{h(+\i), h(-\i)\}$ or $k(+\i)\in \{h(+\i), h(-\i)\}$. 
 The choice of $h$, $k$ implies $s\in \{h(+\i), h(-\i)\}$. 
      Since   
 $d_T(s, k(+\i))=d_T(s, k(-\i))=\pi/2$,    the proposition follows from
  Lemma \ref{6.1lemma3}.

\qed

\subsubsection{When the Minimal Set Contains Singular Points} \label{subsecspecial}

In this section we prove Theorem \ref{gsingular}    when 
$M$ contains   singular points:

\b{Prop}\label{pro6.2}
{Theorem \ref{gsingular}  holds if    there   exists a   special  element  $h_0\in H$ such that 
  $M\cap \p_T {P_{h_0}}$   contains  a  singular point.}

\end{Prop}

Notice Theorem \ref{gsingular}  follows from Propositions    \ref{nonspesin}
     and \ref{pro6.2}.  Also   recall   Lemmas \ref{6.1lemma1},
\ref{6.1lemma2},
\ref{6.1lemma3}   and 
\ref{6.1lemma4}   are still valid.

\b{Le}\label{6.2lemma1}
{If for any special $h\in H$, $M\cap \{h(+\i), h(-\i)\}=\phi$, then 
$H$  either  is   virtually 
free abelian or contains a free group of rank two.}
\end{Le}

\b{proof}
By Lemma \ref{6.1lemma4}  we  may assume    there are two special elements $h,k\in H$
  with 
$\{h(+\i), h(-\i)\}\cap \{k(+\i), k(-\i)\}=\phi$.

Let $m\in M\cap \p_T {P_{h_0}}$ be  a  singular point,  and  
 $g\in H$ be  an arbitrary special element.  Then  
 $d_T(m, g(+\i))\le \pi$, $d_T(m, g(-\i))\le \pi$.   The assumption implies 
$d_T(m, g(+\i))\ge \pi/2$, $d_T(m, g(-\i))\ge \pi/2$. 
If  $d_T(m, g(-\i))=\pi$,  then   Theorem \ref{hyperbolicf} implies 
 $g(+\i)$ is an accumulation point of $\{g^i(m):i\ge 1\}$ and thus lies in $M$,
 contradicting to the assumption.   Therefore $d_T(m, g(-\i))=\pi/2$.
 Similarly $d_T(m, g(+\i))=\pi/2$.

Now $d_T(m, h(+\i))=d_T(m, h(-\i))=d_T(m, k(+\i))=d_T(m, k(-\i))=\pi/2$
  and $\{h(+\i), h(-\i)\}\cap \{k(+\i), k(-\i)\}=\phi$.
The lemma  follows from Lemma \ref{6.1lemma1}.

\end{proof}

\b{Le}\label{6.2lemma2}
{Let $h,k\in H$ be two special elements. If $M\cap \{h(+\i), h(-\i)\}\not=\phi$
   and $M\cap \{k(+\i), k(-\i)\}=\phi$,  then 
$H$  either  is  virtually 
free abelian or contains a free group of rank two.}
\end{Le}
 
\b{proof}
We may assume $h(+\i)\in M$.
Since $M\cap \{k(+\i), k(-\i)\}=\phi$  
the proof of Lemma \ref{6.2lemma1} shows 
 $d_T(h(+\i), k(+\i))=d_T(h(+\i), k(-\i))=\pi/2$. The lemma now follows from 
Lemma \ref{6.1lemma3}.

\end{proof}

\noindent
{\bf{Proof of Proposition \ref{pro6.2}.}}
Suppose $H$  neither   is  virtually 
free abelian nor contains a free group of rank two.
By Lemma \ref{6.1lemma4}     there are two special elements $h,k\in H$
  with 
$\{h(+\i), h(-\i)\}\cap \{k(+\i), k(-\i)\}=\phi$.
Lemmas \ref{6.2lemma1}
    and \ref{6.2lemma2}
  imply that $M\cap \{h(+\i), h(-\i)\}\not=\phi$  and $M\cap \{k(+\i), k(-\i)\}\not=\phi$.

We first claim  $\{h(+\i), h(-\i)\}\subset M$. 
Suppose the claim is false.  After  possibly  replacing $h$ and $k$ with their inverses   
 we   may   assume 
$h(-\i)\in M$, $h(+\i)\notin M$  and $k(-\i)\in M$. 
 We observe $d_T(h(-\i), k(-\i))=\pi/2$ otherwise $h(+\i)$ would be an accumulation point of 
 $\{h^i(k(-\i)):i\ge 1\}$ and thus lie in $M$.  Notice $d_T(h(+\i), k(-\i))\le \pi$.
 If $d_T(h(+\i), k(-\i))=\pi$ 
  and $\sigma$ is a minimal geodesic from $h(+\i)$  to  $k(-\i)$, then $\sigma\cup k(-\i)h(-\i)$ is a
 geodesic from $h(+\i)$  to  $h(-\i)$  with length $1.5\pi$.  It follows  that there is a 
  simple  closed geodesic 
 with length $2.5\pi$ in $\p_T X$ (as $d_T(h(+\i), h(-\i))=\pi$), a contradiction.
If $d_T(h(+\i), k(-\i))=\pi/2$, then we have  $d_T(h(+\i), k(-\i))=d_T(h(-\i), k(-\i))=\pi/2$
   and by Lemma \ref{6.1lemma3}  $H$  either   is  virtually 
free abelian or contains a free group of rank two.   

  Thus   $\{h(+\i), h(-\i)\}\subset M$, and similarly $\{k(+\i), k(-\i)\}\subset M$.
  By assumption we have  $d_T(\xi, \eta)\le \pi$ for any $\xi\in \{h(+\i), h(-\i)\}$ and 
$\eta\in \{k(+\i), k(-\i)\}$.  If $d_T(\xi, \eta)=\pi$ for any $\xi\in \{h(+\i), h(-\i)\}$ and any 
$\eta\in \{k(+\i), k(-\i)\}$, then $<h,k>$ contains a free group of rank two.
  Therefore we may assume there are  $\xi\in \{h(+\i), h(-\i)\}$ and  
$\eta\in \{k(+\i), k(-\i)\}$ with $d_T(\xi, \eta)=\pi/2$.  Without loss of generality we assume 
 $d_T(h(-\i), k(-\i))=\pi/2$.  Then the argument in the preceding paragraph finishes the proof.

\qed




 \addcontentsline{toc}{subsection}{References}

\end{document}